\documentclass[11pt]{article}
\pdfoutput=1
\usepackage{algorithmic}
\usepackage{lscape}
\usepackage{theorem}
\usepackage[english]{babel}
\usepackage{a4}
\usepackage{epsfig}
\usepackage{amsfonts}
\usepackage{amssymb}
\usepackage{latexsym}
\usepackage{url}
\usepackage{amsmath}
\usepackage{clrscode}
\usepackage{comment}
\usepackage[usenames,dvipsnames,svgnames,table]{xcolor}
\usepackage{tikz}
\usetikzlibrary{arrows, positioning,chains,fit,shapes,calc}
\usetikzlibrary{decorations.pathmorphing}
\usepackage{setspace}
\usepackage{multicol}
\usepackage{float}
\usepackage{graphicx}
\usepackage{multirow}
\usepackage{enumitem}
\usepackage{lineno}

\newtheorem{theorem}{Theorem}[section]

\newtheorem{algorithm}[theorem]{Algorithm}
\theorembodyfont{\rmfamily}

\newtheorem{remark}[theorem]{Remark}

\newtheorem{example}[theorem]{Example}
\theorembodyfont{\rmfamily}

\newcommand{\NN}{\ensuremath{\mathbb{N}}}
\newcommand{\ZZ}{\ensuremath{\mathbb{Z}}}
\newcommand{\QQ}{\ensuremath{\mathbb{Q}}}

\newcommand{\CC}{\ensuremath{\mathbb{C}}}

\newcommand{\ii}{\ensuremath{\mathbf{i}}}

\newcommand{\SO}{\ensuremath{\mathop{SO}}}
\newcommand{\Sym}{\ensuremath{\mathrm{S}}}
\newcommand{\Alt}{\ensuremath{\mathrm{A}}}

\definecolor{blue(pigment)}{rgb}{0.2, 0.2, 0.8}
\newcommand{\defn}[1]{{\em \color{blue(pigment)} #1}}

\makeatletter
\newcommand\xleftrightarrow[2][]{\ext@arrow 0099{\longleftrightarrowfill@}{#1}{#2}}
\def\longleftrightarrowfill@{\arrowfill@\leftarrow\relbar\rightarrow}
\makeatother

\tikzset{
  text style/.style={
    sloped, 
    text=black
  }
}

\begin{document}
\title{Solving polynomial systems via homotopy continuation and monodromy}
\newcommand*\samethanks[1][\value{footnote}]{\footnotemark[#1]}
\author{
Timothy Duff \thanks{Research of TD, CH, KL, and AL is supported in part by NSF grants DMS-1151297 and DMS-1719968.} \\ Georgia Tech \and
Cvetelina Hill \samethanks \\ Georgia Tech \and
Anders Jensen\\ Aarhus Universitet \and
Kisun Lee \samethanks \\ Georgia Tech \and
Anton Leykin \samethanks
\\ Georgia Tech \and
Jeff Sommars\\ University of Illinois at Chicago
}
\maketitle

\begin{abstract}
We study methods for finding the solution set of a generic system in a family of polynomial systems with parametric coefficients. 
We present a framework for describing monodromy based solvers in terms of decorated graphs. 
Under the theoretical assumption that monodromy actions are generated uniformly, 
we show that the expected number of homotopy paths tracked by an algorithm following this framework is linear in the number of solutions. We demonstrate that our software implementation is competitive with the existing state-of-the-art methods implemented in other software packages.

\end{abstract}
\tableofcontents

\section{Introduction}
Homotopy continuation has become a standard technique to find approximations of solutions of polynomial systems. There is an early popular text on the subject and its applications by Morgan~\cite{Morgan87}. This technique is the backbone of Numerical Algebraic Geometry, the area which classically addresses the questions of complex algebraic geometry through algorithms that  employ numerical approximate computation. The chapter by Sommese, Verschelde, and Wampler~\cite[\S 8]{SVW9} is the earliest introduction and the book by Sommese and Wampler~\cite{Sommese-Wampler-book-05} is the primary reference in the area.       

Families of polynomial systems with parametric coefficients play one of the central roles. Most homotopy continuation techniques could be viewed as going from a generic system in the family to a particular one. This process is commonly referred to as \defn{degeneration}. Going in the reverse direction, it may be called \defn{deformation}, \defn{undegeneration}, or \defn{regeneration} depending on the literature. Knowing the solutions of a generic system one can use \defn{coefficient-parameter homotopy}~\cite[\S 7]{Sommese-Wampler-book-05} to get to the solution of a particular one.

The main problem that we address here is how to solve a generic system in a family of systems  
$$F_p = (f_p^{(1)},\ldots,f_p^{(N)}) = 0, \quad f_p^{(i)}\in \CC[p][x],\ i=1,\ldots,N,$$ 
with finitely many parameters $p$ and $n$ variables $x$.
In the main body of the paper we restrict our attention to \defn{linear parametric} families of systems, defined as systems with affine linear parametric coefficients, such that for a generic $p$ we have a nonempty finite set of solutions $x$ to $F_p(x)=0$. This implies $N\geq n$.  The number of parameters is arbitrary, but we require that for a generic $x$ there exists $p$ with $F_p(x)=0$.
These restrictions are made for the sake of simplicity. We explain what modifications are needed to apply our approach in more general settings in~\S\ref{sec:generalizations}.

Linear parametric systems form a large class that includes \defn{sparse polynomial systems}. These are square ($n=N$) systems with a fixed monomial support for each equation and a distinct parameter for the coefficient of each monomial. \defn{Polyhedral homotopy} methods for solving sparse systems stem from the BKK (Bernstein, Khovanskii, Kouchnirenko) bound on the number of solutions~\cite{Bernstein:BKK}; the early work on algorithm development was done  in~\cite{HuberSturmfels:PolyhedralHomotopies,Verschelde-Verlinden-Cools}. Polyhedral homotopies provide an optimal solution to sparse systems in the sense that they are designed to follow exactly as many paths as the number of solutions of a generic system (the BKK bound). 

The method that we propose is clearly not optimal in the above sense. The expected number of homotopy paths followed can be larger than the number of solutions, though \emph{not significantly larger}. We also use linear segment homotopies that are \emph{significantly simpler} and less expensive to follow in practice. Our current implementation shows it is competitive with the state-of-the-art implementations of polyhedral homotopies in PHCpack~\cite{V99} and HOM4PS2~\cite{lee2008hom4ps} for solving sparse systems. In a setting more general than sparse, we demonstrate examples of linear parametric systems for which our implementation exceeds the capabilities of the existing sparse system solvers and blackbox solvers based on other ideas.      

The idea of using the monodromy action induced by the fundamental group of the regular locus of the parameter space has been successfully employed throughout Numerical Algebraic Geometry. One of the main tools in the area, \defn{numerical irreducible decomposition}, can be efficiently implemented using the \defn{monodromy breakup} algorithm, which first appeared in~\cite{SVW2001:monodromy}. One parallel incarnation of the monodromy breakup algorithm is described in~\cite{Ley-Ver-new-monodromy-05}. In fact, the main idea in that work is close in spirit to what we propose in this article. The idea to use monodromy to find solutions drives numerical implicitization~\cite{Chen-Kileel:NumericalImplicitization} and appears in other works such as~\cite{del2015critical}. Computing monodromy groups numerically, as in~\cite{Leykin-Sottile:HoG} and \cite{hauenstein2016numerical}, requires more computation than just finding solutions. One can approach this computation with the same methodology as we propose; see (\ref{sec:Galois}) of \S\ref{sec:generalizations}. 

Our main contribution is a new framework to describe algorithms for solving polynomial systems using monodromy; we call it the \defn{Monodromy Solver} (MS) framework. We analyze the complexity of our main algorithm both theoretically assuming a certain statistical model and experimentally on families of examples. The analysis gives us grounds to say that the expected number of paths tracked by our method is linear, with a small coefficient, as the number of solutions grows.
Our method and its implementation not only provide a new general tool for solving polynomial systems, but also can solve some problems out of reach for other existing software.

The structure of the paper is as follows. We give a brief overview of the MS method intermingled with some necessary preliminaries in~\S\ref{sec:overview}. 
An algorithm following the MS framework depends on a choice of strategy, with several possibilities outlined in~\S\ref{sec:strategies}. 
Statistical analysis of the method is the topic of~\S\ref{sec:random_loops}. 
The implementation is discussed in~\S\ref{sec:implementation} together with the side topic of \defn{certification} of the solution set. 
The results of our experiments on selected example families highlighting various practical computational aspects are in~\S\ref{sec:experiments}. The reader may also want look at examples of systems in~\S\ref{subsec:sparseexperiments}~and~\S\ref{subsec:CRN} before reading some earlier sections. 
Possible generalizations of the MS technique and the future directions to explore are presented in~\S\ref{sec:generalizations}.

\section{Background and framework overview}\label{sec:overview}
Let $m, n\in\NN$. We consider the complex linear space of square systems $F_p$, $p\in \CC^m$, 
where the monomial support of $f_p^{(1)},\dots,f_p^{(n)}$ in the variables $x=(x_1,\dots,x_n)$ is fixed and the coefficients vary.
By a \defn{base space} $B$ we mean a parametrized linear variety of systems. We think of it as the image of an affine linear map $\varphi: p\mapsto F_p$ from a parameter space $\CC^m$ with coordinates $p=(p_1,\ldots,p_m)$ to the space of systems.

We assume the structure of our family is such that the projection $\pi$ from the \defn{solution variety} 
\[
V  = \{(F_p,x) \in B\times\CC^n \mid F_p(x)=0\}
\]
to $B$ gives us a \defn{branched covering}, i.e., the fiber $\pi^{-1}(F_p)$ is finite of the same cardinality for a generic $p$.
The discriminant variety $D$ in this context is the subset of the systems in the base space with nongeneric fibers; it is also known as the \defn{branch locus} of $\pi$.

The \defn{fundamental group} $\pi_1(B\setminus D)$ --- note that $\pi_1$ is a usual topological notation that is not related to the map $\pi$ above --- as a set consists of loops, i.e., paths in $B\setminus D$ starting and finishing at a fixed $p\in B\setminus D$ considered up to homotopy equivalence. The definition, more details to which one can find in  \S\ref{sec:monodromy}, does not depend on the point $p$, since $B\setminus D$ is connected. Each loop induces a permutation of the fiber $\pi^{-1}(F_{p})$, which is referred to as a \defn{monodromy} action. 

Our goal is to find the fiber of one generic system in our family.
Our method is to find one pair  $(p_0,x_0)\in V$ 
and use the monodromy action on the fiber $\pi^{-1}(F_{p_0})$ to find its points. We assume that this action is transitive, which is the case if and only if the solution variety $V$ is irreducible. If $V$ happens to be reducible, we replace $V$ with its unique dominant irreducible component as explained in Remark~\ref{rem:dominant-component}.

\subsection{Monodromy}\label{sec:monodromy}

We briefly review the basic facts concerning monodromy groups of branched coverings. With notation as before, fix a system $F_p\in B \setminus D$ and consider a loop $\tau $ without branch points based at $F_p$; that is, a continuous path 
\[
\tau : [0,1] \rightarrow B \setminus D
\] 
such that $\tau (0) = \tau (1) = F_p.$ Suppose we are also given a
point $x_i$ in the fiber $\pi^{-1} (F_p)$ with $d$ points $x_1, x_2, \ldots , x_d.$ Since $\pi $ is a covering map, the pair $(\tau, x_i)$ corresponds to a unique \defn{lifting} $\widetilde{\tau_i},$ a path
\[
\widetilde{\tau_i} : [0,1] \rightarrow V
\] 
such that
$\widetilde{\tau_i} (0) = x_i$ and
$\widetilde{\tau_i} (1) = x_j$ for some $1\le j \le d.$
Note that the reversal of $\tau$ and $x_j$ lift to a reversal of $\widetilde{\tau_i}$.
Thus, the loop $\tau $ induces a permutation of the set $\pi^{-1} (F_p).$ We have a group homomorphism
\[
\varphi : \pi_1 (B\setminus D, F_p) \rightarrow \Sym_d
\]
whose domain is the usual fundamental group of $B \setminus D$ based at $F_p$. The image of $\varphi $ is the \defn{monodromy group} associated to $\pi^{-1} (F_p).$ The monodromy group acts on the fiber $\pi^{-1} (F_p)$ by permuting the solutions of $F_p$.

\begin{remark}
A reader familiar with the notion of a \defn{monodromy loop} in the discussion of~\cite[\S
15.4]{Sommese-Wampler-book-05} may think of this keyword referring to
a representative of an element of the fundamental group together with
its liftings to the solution variety and the induced action on the fiber.
For the purposes of this article we need to be clear about the
ingredients bundled in this term.
\end{remark}

We have not used any algebraic properties so far. The construction of
the monodromy group above holds for an arbitrary covering with finitely many sheets. The monodromy group is a transitive subgroup of $\Sym_d$ whenever the total space is connected.
In our setting, since we are working over $\CC$, this occurs precisely when the solution variety is irreducible.
\begin{remark}\label{rem:dominant-component}  For a linear family, we can show that there is at most one irreducible component of the solution variety $V$ for which the restriction of the projection $(F_p,x)\mapsto x$ is dominant (that is, its image is dense). We call such component the \defn{dominant component}.  
Indeed, let $U$ be the locus of points $(F_p,x)\in \pi^{-1}(B\setminus D)$ such that 
\begin{itemize}
\itemsep=0mm
\item the restriction of the $x$-projection map is locally surjective, and
\item the solution to the linear system of equations $F_p(x)=0$ in $p$ has the generic dimension. 
\end{itemize}
Being locally surjective could be interpreted either in the sense of Zariski topology or as inducing surjection on the tangent spaces.
Then either $U$ is empty or $\overline{U}$ is the dominant component we need, since it is a
vector bundle over an irreducible variety, and is hence irreducible. 
\end{remark}

In the rest of the paper, when we say \emph{solution variety},  we mean the \emph{dominant component of the solution variety}.
In particular, for sparse systems restricting the attention to the dominant component translates into looking for solutions only in the torus $(\CC^*)^n$.

\subsection{Homotopy continuation}\label{subsec:homotopy}
Given two points $F_{p_1}$ and $F_{p_2}$ in the base space $B$, we may form the family of systems 
$$
H(t) = (1-t)F_{p_1} + t F_{p_2}\,,\quad t\in[0,1],
$$ 
known as the \defn{linear segment homotopy} between the two systems. 
If $p_1$ and $p_2$ are sufficiently generic, for each $t\in[0,1]$ we have $H(t)$ outside the real codimension 2 set $D$. Consequently, each system $H(t)$ has a finite and equal number of solutions. This homotopy is a path in $B$; a lifting of this path in the solution variety $V$ is  called a \defn{homotopy path}. The homotopy paths of $H(t)$ establish a one-to-one correspondence between the fibers $\pi^{-1}(F_{p_1})$ and $\pi^{-1}(F_{p_2})$.

\begin{remark}\label{rem:gamma-trick}
Note that $\gamma F_p$ for $\gamma\in\CC\setminus\{0\}$  has the same solutions as $F_p$.
Let us scale both ends of the homotopy by taking a homotopy between $\gamma_1F_{p_1}$ and   $\gamma_2F_{p_2}$ for generic $\gamma_1$ and $\gamma_2$.
If the coefficients of $F_p$ are homogeneous in $p$ then 
$$
H'(t) = (1-t)\gamma_1F_{p_1} + t \gamma_2F_{p_2}  = F_{(1-t)\gamma_1p_1+t\gamma_2p_2}\,, \quad t\in[0,1],
$$
is a homotopy matching solutions $\pi^{-1}(F_{p_1})$ and  $\pi^{-1}(F_{p_2})$ where the matching is potentially different from that given by $H(t)$.
 Similarly, for an affine linear family, $F_p = F'_p + C$ where $F'_p$ is homogeneous in $p$ and $C$ is a constant system, we have 
$$
H'(t) =  (1-t)\gamma_1F_{p_1} + t \gamma_2F_{p_2}  = F'_{(1-t)\gamma_1p_1+t\gamma_2p_2}+((1-t)\gamma_1+t\gamma_2)C.
$$
We ignore the fact that $H'(t)$ may go outside $B$ for $t\in(0,1)$, since its rescaling,
\begin{align*}
H''(t) &= \frac{1}{(1-t)\gamma_1+t\gamma_2}H'(t)\\
&=F'_{\frac{(1-t)\gamma_1p_1+t\gamma_2p_2}{(1-t)\gamma_1+t\gamma_2}}+C = F_{\frac{(1-t)\gamma_1p_1+t\gamma_2p_2}{(1-t)\gamma_1+t\gamma_2}}\,, \quad t\in[0,1],
\end{align*}
does not leave $B$ and clearly has the same homotopy paths. Note that $H''(t)$ is well defined as $(1-t)\gamma_1+t\gamma_2\neq 0$ for all $t\in[0,1]$ for generic $\gamma_1$ and $\gamma_2$.
\end{remark}

One may use methods of \defn{numerical homotopy continuation}, described, for instance, in~\cite[\S 2.3]{Sommese-Wampler-book-05}, to track the solutions as $t$ changes from $0$ to $1$. In some situations the path in $B$ may pass close to the branch locus $D$ and numerical issues must be considered.

\begin{remark}\label{rem:parameter-segment-homotopy}
If the family $F_p$ is nonlinear in the parameters $p$, one has to take the \defn{parameter linear segment homotopy} in the parameter space, i.e., $H(t) = F_{(1-t)p_1+tp_2}$, $t\in[0,1]$. This does not change the overall construction; however, the freedom to replace the systems $F_{p_1}$ and $F_{p_2}$ at the ends of the homotopy with their scalar multiples as in Remark~\ref{rem:gamma-trick} is lost. 
\end{remark}

\subsection{Graph of homotopies: main ideas}
\label{sec:mainIdeas}
Some readers may find it helpful to use the examples of \S\ref{sec:graph-of-homotopies-examples} for graphical intuition as we introduce notation and definitions below.  
\smallskip

To organize the discovery of new solutions we represent the set of homotopies by a finite undirected graph $G$. Let $E(G)$ and $V(G)$ denote the edge and vertex set of $G$, respectively. Any vertex $v$ in $V(G)$ is associated to a point $F_p$ in the base space. An edge $e$ in $E(G)$ connecting $v_1$ and $v_2$ in $V(G)$ is decorated with two complex numbers, $\gamma_1$ and $\gamma_2$, and represents the linear homotopy connecting $\gamma_1F_{p_1}$ and $\gamma_2 F_{p_2}$ along a line segment (Remark~\ref{rem:gamma-trick}).
We assume that both $p_i$ and $\gamma_i$ are chosen so that the segments do not intersect the branch locus. Choosing these at random (see~\S\ref{sec:randomization} for a possible choice of distribution) satisfies the assumption, since the exceptional set of choices where such intersections happen is contained in a real Zariski closed set, see \cite[Lemma 7.1.3]{Sommese-Wampler-book-05}.

We allow multiple edges between two distinct vertices but no loops, since the latter induce trivial homotopies.
For a graph $G$ to be potentially useful in a monodromy computation, it must contain a cycle. Some 
of the general ideas behind the structure of a graph $G$ are listed below. 
\begin{itemize}
\item For each vertex $v_i$, we maintain a subset of \defn{known} points $Q_i\subset\pi^{-1}(F_{p_i})$.
\item For each edge $e$ between $v_i$ and $v_j$, we record the two complex numbers $\gamma_1$ and $\gamma_2$ and we store the known partial correspondences $C_e\subset\pi^{-1}(F_{p_i})\times\pi^{-1}(F_{p_j})$ 
between known points $Q_i$ and $Q_j$.  
\item At each iteration, we pick an edge and direction, track the corresponding homotopy starting with yet unmatched points, and update known points and correspondences between them.
\item We may obtain the initial ``knowledge'' as a \defn{seed pair} $(p_0,x_0)$ by picking $x_0\in\CC^n$ at random and choosing $p_0$ to be a generic solution of the linear system $F_{p}(x_0)=0$.
\end{itemize}
We list basic operations that result in transition between one state of our algorithm captured by $G$, $Q_i$ for $v_i\in V(G)$, and $C_e$ for $e\in E(G)$ to another.  
\begin{enumerate}
\item \label{operation:track}For an edge $e = v_i \xleftrightarrow{(\gamma_1,\gamma_2)} v_j$, consider the homotopy 
$$H^{(e)} = (1-t)\gamma_1F_{p_i} + t \gamma_2F_{p_j}$$ 
where $(\gamma_1,\gamma_2)\in\CC^2$ is the label of $e$. 
    \begin{itemize} 
         \item Take start points $S_i$ to be a subset of the set of known points $Q_i$ that do not have an established correspondence with points in $Q_j$.  
         \item Track $S_i$ along $H^{(e)}$ for $t\in[0,1]$ to get $S_j \subset \pi^{-1}(F_{p_j})$.
         \item Extend the known points for $v_j$, that is, $Q_j := Q_j \cup S_j$ and record the newly established correspondences.   
    \end{itemize}
\item \label{operation:add-vertex} Add a new vertex corresponding to $F_{p}$ for a generic $p \in B\setminus D$.
\item \label{operation:add-edge} Add a new edge $e = v_i \xleftrightarrow{(\gamma_1,\gamma_2)} v_j$ between two existing vertices decorated with generic $\gamma_1,\gamma_2\in\CC$.
\end{enumerate}

At this point a reader who is ready to see a more formal algorithm based on these ideas may skip to Algorithm~\ref{alg:sgs}.  

\subsection{Graph of homotopies: examples}\label{sec:graph-of-homotopies-examples}
We demonstrate the idea of graphs of homotopies, the core idea of the MS framework, by giving two examples.

\begin{example}
Figure~\ref{fig:BigPicture} shows a graph $G$ with 2 vertices and 3 edges embedded in the base space $B$ with paths partially lifted to the solution variety, which is a covering space with 3 sheets. 
The two fibers $\{x_1,x_2,x_3\}$ and $\{y_1,y_2,y_3\}$ are connected
by 3 partial correspondences induced by the liftings of three
egde-paths. 

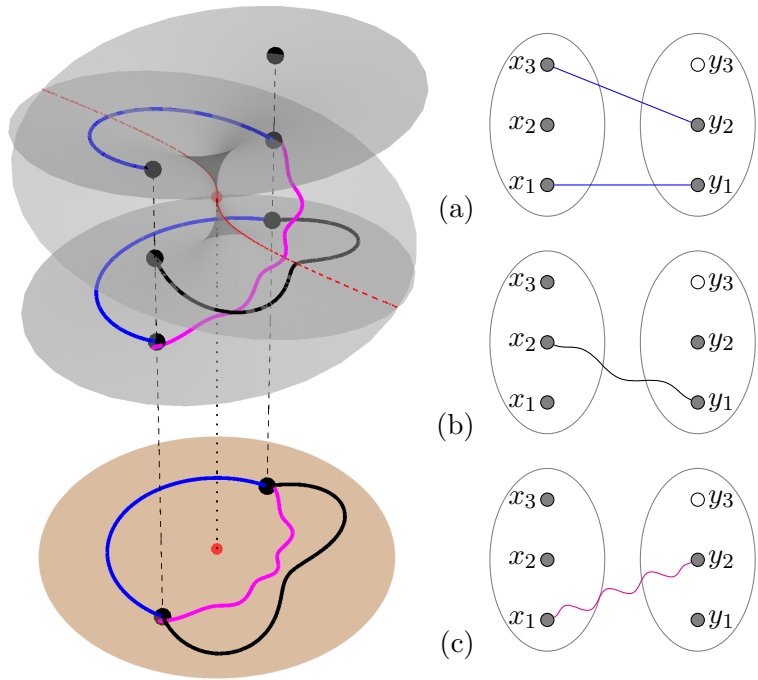
\begin{figure}[H]
\vspace{-0.5cm}
\begin{multicols}{2}
        \null \vfill
	\centering
	\includegraphics[width=8cm]{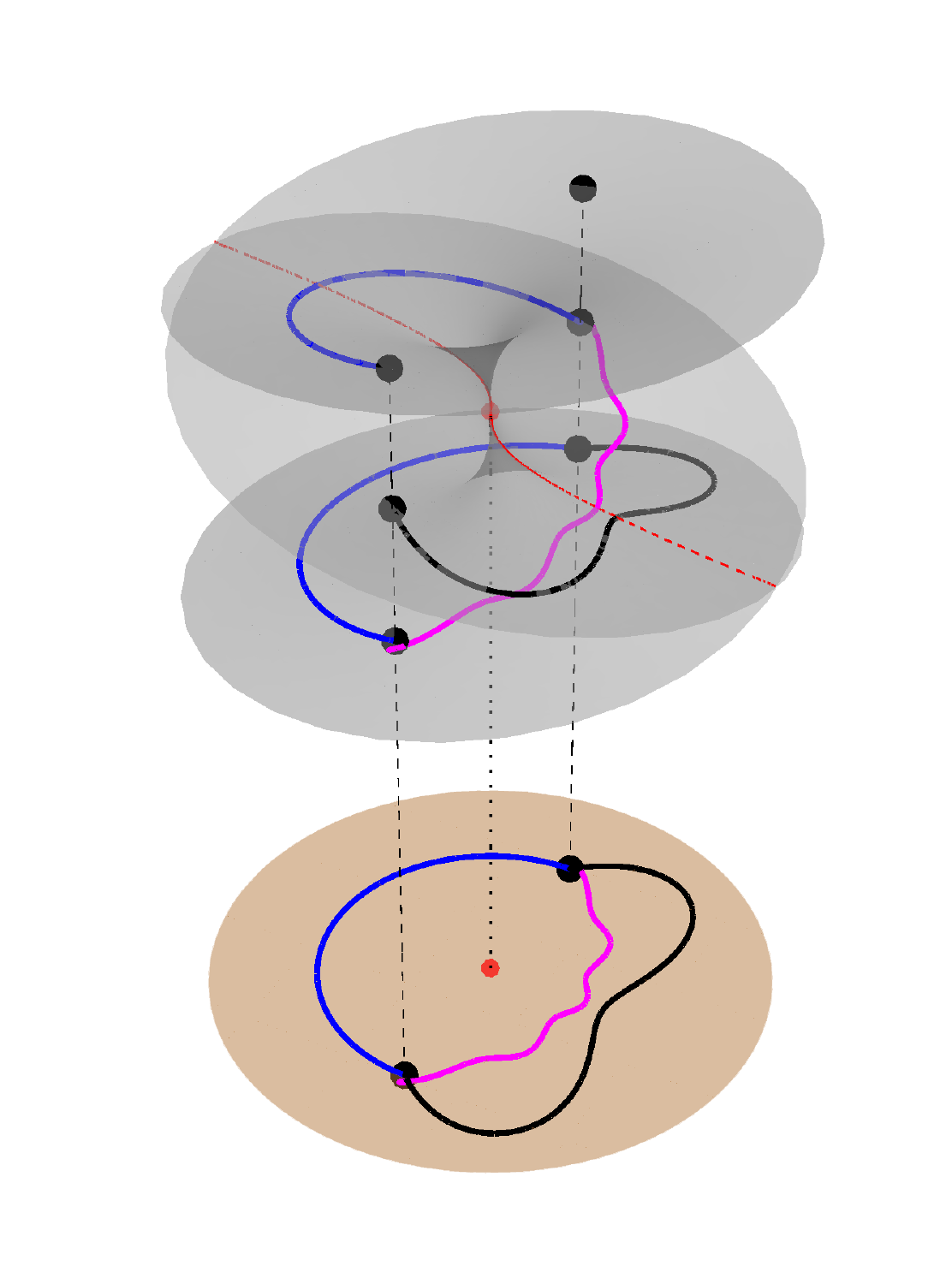}
        \begin{enumerate}
        \vfill \null 
        \columnbreak
        \null \vfill
        \item[(a)]
            \begin{tikzpicture}[auto, every node/.style={circle, draw,
               inner sep=0pt,minimum size=5pt, node distance=8mm}, every
             fit/.style={ellipse,draw,inner sep=-1pt,text
               centered, text width=1.15cm}]
             \node (a1) [label=left:$x_3$,fill=gray] {};
             \node (a2) [below of=a1,label=left:$x_2$, fill=gray] {};
             \node (a3) [below of=a2,label=left:$x_1$, fill=gray] {};
             \node[gray, fit=(a1)(a3)] {};
             
             \node (b1) [xshift=2cm,label=right:$y_3$] {};
             \node (b2) [below of=b1,label=right:$y_2$, fill=gray] {};
             \node (b3) [below of=b2,label=right:$y_1$, fill=gray] {};
             \node[gray, fit=(b1)(b3)] {};

             \draw (a3) -- (b3) [blue];
             \draw (a1) -- (b2) [blue];
            \end{tikzpicture}

        \item[(b)]
            \begin{tikzpicture}[auto, every node/.style={circle, draw,
               inner sep=0pt,minimum size=5pt, node distance=8mm}, every
             fit/.style={ellipse,draw,inner sep=-1pt,text
               centered, text width=1.15cm}]
             \node (a1) [label=left:$x_3$,fill=gray] {};
             \node (a2) [below of=a1,label=left:$x_2$, fill=gray] {};
             \node (a3) [below of=a2,label=left:$x_1$, fill=gray] {};
             \node[gray, fit=(a1)(a3)] {};
             
             \node (b1) [xshift=2cm,label=right:$y_3$] {};
             \node (b2) [below of=b1,label=right:$y_2$, fill=gray] {};
             \node (b3) [below of=b2,label=right:$y_1$, fill=gray] {};
             \node[gray, fit=(b1)(b3)] {};

             \draw [black,decorate, decoration=snake, segment length=12mm] (a2) -- (b3);
            \end{tikzpicture}

        \item[(c)]
            \begin{tikzpicture}[auto, every node/.style={circle, draw,
               inner sep=0pt,minimum size=5pt, node distance=8mm}, every
             fit/.style={ellipse,draw,inner sep=-1pt,text
               centered, text width=1.15cm}]
             \node (a1) [label=left:$x_3$,fill=gray] {};
             \node (a2) [below of=a1,label=left:$x_2$, fill=gray] {};
             \node (a3) [below of=a2,label=left:$x_1$, fill=gray] {};
             \node[gray, fit=(a1)(a3)] {};
             
             \node (b1) [xshift=2cm,label=right:$y_3$] {};
             \node (b2) [below of=b1,label=right:$y_2$, fill=gray] {};
             \node (b3) [below of=b2,label=right:$y_1$, fill=gray] {};
             \node [gray, fit=(b1)(b3)] {};

             \draw[magenta, decorate, decoration=snake, segment length=6mm] (a3) -- (b2);
           \end{tikzpicture}
        \end{enumerate}
        \vfill \null
\end{multicols}
\vspace{-1.5cm}
	\caption{Selected liftings of 3 edges connecting the fibers of 2 vertices and induced correspondences.}
	\label{fig:BigPicture}
\end{figure}

Note that several aspects in this illustration are fictional. There
is only one branch point in the actual complex base space $B$
that we would like the reader to imagine. The visible
self-intersections of the solution variety $V$ are an artifact of
drawing the picture in the real space. Also, in practice we use
homotopy paths as simple as possible, 
however, here the paths are more involved for the purpose of distinguishing them in print.   

An algorithm that we envision may hypothetically take the following steps:
\begin{enumerate}[label={(\arabic*)}]
\itemsep=0mm
\item \emph{seed} the first fiber with $x_1$;
\item use a lifting of edge $e_a$ to get $y_1$ from $x_1$;
\item use a lifting of edge $e_b$ to get $x_2$ from $y_1$;
\item use a lifting of edge $e_c$ to get $y_2$ from $x_1$;
\item use a lifting of edge $e_a$ to get $x_3$ from $y_2$.
\end{enumerate}
Note that it is \emph{not} necessary to complete the correspondences
(a), (b), and (c). Doing so would require tracking 9 continuation
paths, while the hypothetical run above uses only 4 paths to find a fiber.
\end{example}

\begin{example}\label{example:two-correspondences}
Figure~\ref{fig:correspondences} illustrates two partial
correspondences associated to two edges $e_a$ and $e_b$, both
connecting two vertices $v_1$ and $v_2$ in $V(G)$. Each vertex $v_i$
stores the array of known points $Q_i$, which are depicted in solid. 
Both correspondences in the picture are subsets of a perfect matching,
a one-to-one correspondence established by a homotopy associated to the edge. 
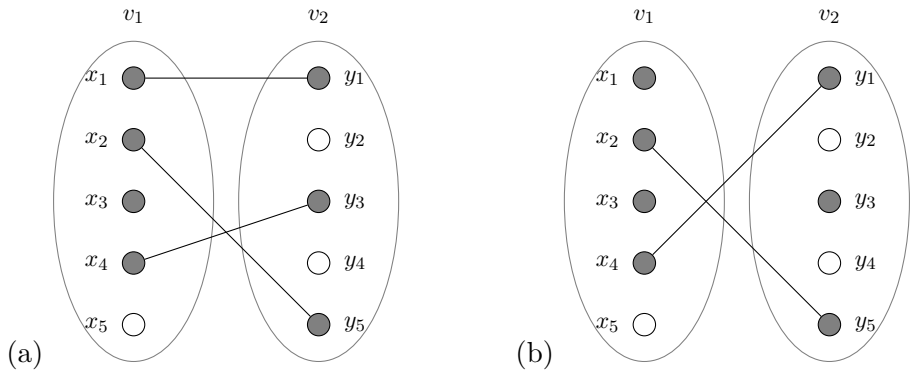
\begin{figure}[H]
\begin{multicols}{2}
  \centering
  (a) \begin{tikzpicture}[auto, every node/.style={circle,draw,scale=0.82},
    every fit/.style={ellipse,draw,inner sep=1.2pt,text width=1.75cm}, scale=0.82]
    \node (a0) [label=left:$x_1$, fill=gray]	{};
    \node (a1) [below of=a0,label=left:$x_2$, fill=gray] {};
    \node (a2) [below of=a1,label=left:$x_3$, fill=gray] {};
    \node (a3) [below of=a2,label=left:$x_4$, fill=gray] {};
    \node (a4) [below of=a3,label=left:$x_5$] {};
    \node [gray,fit=(a0) (a4),label=above:$v_1$] {};
    
    \node (b0) [xshift=3cm,label=right:$y_1$, fill=gray]	{};
    \node (b1) [below of=b0,label=right:$y_2$] {};
    \node (b2) [below of=b1,label=right:$y_3$, fill=gray] {};
    \node (b3) [below of=b2,label=right:$y_4$] {};
    \node (b4) [below of=b3,label=right:$y_5$, fill=gray] {};
    \node [gray,fit=(b0) (b4),label=above:$v_2$] {};
                
    \draw (a0) -- (b0);
    \draw (a3) -- (b2);
    \draw (a1) -- (b4);
  \end{tikzpicture}

  (b) \begin{tikzpicture}[auto, every node/.style={circle,draw, scale=0.82},
    every fit/.style={ellipse,draw,inner sep=1.2pt,text width=1.75cm}, scale=0.82]
    \node (a0) [label=left:$x_1$, fill=gray]	{};
    \node (a1) [below of=a0,label=left:$x_2$, fill=gray] {};
    \node (a2) [below of=a1,label=left:$x_3$, fill=gray] {};
    \node (a3) [below of=a2,label=left:$x_4$, fill=gray] {};
    \node (a4) [below of=a3,label=left:$x_5$] {};
    \node [gray,fit=(a0) (a4),label=above:$v_1$] {};
    
    \node (b0) [xshift=3cm,label=right:$y_1$, fill=gray]	{};
    \node (b1) [below of=b0,label=right:$y_2$] {};
    \node (b2) [below of=b1,label=right:$y_3$, fill=gray] {};
    \node (b3) [below of=b2,label=right:$y_4$] {};
    \node (b4) [below of=b3,label=right:$y_5$, fill=gray] {};
    \node [gray,fit=(b0) (b4),label=above:$v_2$] {};
                
    \draw (a3) -- (b0);
    \draw (a1) -- (b4);
  \end{tikzpicture}

\end{multicols}
		\caption{Two partial correspondences induced by edges $e_a$ and $e_b$ for the fibers of the covering map of degree $d=5$ in Example~\ref{example:two-correspondences}.}
		\label{fig:correspondences}
\end{figure}
Note that taking the set of start points $S_1=\{x_3\}$ and following the homotopy $H^{(e_a)}$ from left to right is \emph{guaranteed} to discover a new point in the second fiber. On the other hand, it is impossible to obtain new knowledge by tracking $H^{(e_a)}$ from right to left. 
Homotopy $H^{(e_b)}$ has a \emph{potential} to discover new points if
tracked in either direction. We can choose $S_1=\{x_1,x_3\}$ as the
start points for one direction and $S_2=\{y_3\}$ for the other. In
this scenario, following the homotopy from left to right is guaranteed
to produce at least one new point, while going the other way may
either deliver a new point or just augment the correspondences between
the already known points. If the correspondences in (a) and (b) are
completed to one-to-one correspondences of the fibers, taking the
homotopy induced by the edge $e_a$ from left to right followed by the
homotopy induced by edge $e_b$ from right to left would produce a
permutation. However, the group generated by this permutation has to
stabilize $\{x_2\}$, therefore, it would not act transitively on the
fiber of $v_1$. One could also imagine a completion such that the
given edges would not be sufficient to discover $x_5$ and~$y_4$. 
\end{example}

In our algorithm, we record and use correspondences; however, they are viewed as a secondary kind of knowledge. In particular,  in~\S\ref{sec:edgeselection} we develop heuristics driven by edge potential functions which look to maximize the number of newly discovered solutions, in other words, to extend the primary knowledge in some greedy way.

\section{Algorithms and strategies}\label{sec:strategies}
The operations listed in \S\ref{sec:mainIdeas} give a great deal
of freedom in the discovery of solutions. However, not all strategies
for applying these operations are equally efficient. We distinguish
between \defn{static} strategies, where the graph is fixed throughout
the discovery process (only basic operation~\ref{operation:track} of
\S\ref{sec:mainIdeas} is used) and \defn{dynamic} strategies, where
vertices and edges may be added
(operations~\ref{operation:add-vertex} and~\ref{operation:add-edge}).

\subsection{A naive dynamic strategy}\label{subsec:naive}
To visualize this strategy in our framework jump ahead and
to the \verb|flower| graph in Figure~\ref{fig:rFS}. Start with
the seed solution at the vertex $v_0$ and proceed creating loops as
petals in this graph: e.g., 
use basic
operations~\ref{operation:add-vertex} and~\ref{operation:add-edge} to
create $v_1$ and two edges between $v_0$ and $v_1$, track the
known solutions at $v_0$ along the new petal to potentially find new
solutions at $v_0$, then ``forget'' the petal and create an
entirely new one in the next iteration.   

This strategy populates the fiber $\pi^{-1}(F_{p_0})$, but how fast?
Assume the permutation induced by a petal permutation on $\pi^{-1}(F_{p_1})$
is uniformly distributed. Then for the first petal the probability of
finding a new solution is $(d-1)/d$ where
$d=|\pi^{-1}(F_{p_1})|$. This probability is close to $1$ when $d$ is large, however for the other petals the probability of
arriving at anything new at the end of one tracked path decreases as the known solution set
grows.
 
Finding the expected number of iterations (petals) to discover the
entire fiber is equivalent to solving the \defn{coupon collector's
  problem}. The number of iterations is $d\,\ell(d)$ where
$\ell(d):={1\over 1}+{1\over 2}+\cdots+{1\over d}$. The values of
$\ell(d)$ can be regarded as lower and upper sums for two integrals of
the function $x\mapsto x^{-1}$, leading to the bounds
$\ln(d+1)\leq \ell(d)\leq \ln(d)+1$.  
Simultaneously tracking all known points along a petal gives a better complexity, since different paths cannot lead to the same solution.

We remark that the existing implementations of numerical irreducible decomposition in Bertini~\cite{Bertini-book}, PHCpack~\cite{V99}, and {\tt NumericalAlgebraicGeometry} for Macaulay2~\cite{Leykin:NAG4M2} that use monodromy are driven by a version of the naive dynamic strategy. 

\subsection{Static graph strategies}
It turns out to be an advantage to reuse the edges of the graph. In a static strategy the graph is fixed and we discover solutions according to the following algorithm.
\begin{algorithm}[Static graph strategy] \label{alg:sgs} Let the base space be given by a map $\varphi:p\mapsto F_p$.

\noindent
$(j,Q_j) = {\tt monodromySolve}(G,Q',{\tt stop})$
\begin{algorithmic}[0]
\renewcommand{\algorithmicrequire}{\textbf{Input:}}
\renewcommand{\algorithmicensure}{\textbf{Output:}}
\REQUIRE
~\\
\begin{itemize}
\itemsep=0mm
\item A graph $G$ with vertices decorated with $p_i$'s and edges decorated with pairs $(\gamma_1,\gamma_2)\in\CC^2$.
\item Subsets $Q'_i\subset\pi^{-1}(\varphi(p_i))$ for $i\in 1,\dots,|V(G)|$, not all empty.
\item A stopping criterion ${\tt stop}$.
\end{itemize}
\ENSURE A vertex $j$ in $G$ and a subset $Q_j$ of the fiber $\pi^{-1}(F_{p_j})$ with the property that $Q_j$ cannot be extended by tracking homotopy paths represented by $G$.

\smallskip\hrule\smallskip

\STATE $Q_i:=Q'_i$ for $i\in 1,\dots,|V(G)|$.
\WHILE{there exists an edge $e=(j,k)$ in $G$ such that $Q_j$ has points not yet tracked with $H^{(e)}$}
 \STATE Choose such an edge $e=(j,k)$.
 \STATE Let $S\subset Q_j$ be a nonempty subset of the set of points not yet tracked with $H^{(e)}$. \label{line:start-points}
 \STATE Track the points $S$ with $H^{(e)}$ to obtain elements $T\subset \pi^{-1}(\varphi(p_k)) \setminus Q_k$.
 \STATE Let $Q_k:=Q_k\cup T$.
 \IF{the criterion ${\tt stop}$ is satisfied (e.g., $|Q_k|$ equals a known solution count)} 
  \RETURN $(k,Q_k)$
 \ENDIF 
\ENDWHILE
\STATE Choose some vertex $j$ and return $(j,Q_j)$.
\smallskip\hrule\smallskip
\end{algorithmic}
\end{algorithm}

The algorithm can be specialized in several ways. We may
\begin{itemize}
\item choose the graph $G$,
\item specify a stopping criterion ${\tt stop}$,
\item choose a strategy for picking the edge $e=(j,k)$.
\end{itemize}
We address the first choice in \S\ref{sec:layouts} by listing several graph layouts that can be used. Stopping criteria are discussed in \S\ref{sec:solutioncount} and \S\ref{sec:no-solution-count}, while strategies for selecting an edge are discussed in \S\ref{sec:edgeselection}.

\begin{remark}\label{rem:complexity}
 We notice that if the stopping criterion is never satisfied, the
 number of paths being tracked by Algorithm~\ref{alg:sgs} is at most
 $d|E(G)|$, where $d$ is the number of solutions of a generic system.
\end{remark}

\subsubsection{Two static graph layouts}
\label{sec:layouts}
\begin{figure}[H]
	\centering
\begin{minipage}{.4\textwidth}
		\begin{tikzpicture}[auto, thick,main
                  node/.style={circle,draw,font=\sffamily\Large\bfseries,scale=0.6},
                  scale=0.6]
		\node[main node] (0) at (0, 0) {$v_0$};
		\node[main node] (1)  at (3.5,2.3) {$v_1$};
		\node[main node] (2) at (3.5,-2.3) {$v_2$};
		\node[main node] (3) at (-3.5, 2.3) {$v_3$};
		\node[main node] (4) at (-3.5, -2.3) {$v_4$};

		\path[every node/.style={font=\sffamily\small}]
		(0) edge [bend right] node [right] {} node [above=+.2, text style] 
			{
} (1) 
		        edge [bend right] node[left] {} node [above=+.2, text style] 
		        {
} (2)
		        edge [bend left] node[left] {} node [above=+.2, text style] 
		        {
} (3)
		        edge [bend left] node[left] {} node [above=+.2, text style] 
		        {
} (4)
		(1) edge [bend right] node[left] {} node [below=-.05 , text style] 
			{
} (0)
		(2) edge [bend right] node[left] {} node [below=-.05, text style]
			{
} (0)
		(3) edge [bend left] node[left] {} node [below=-.05, text style]
			{
} (0)
		(4) edge [bend left] node[left] {} node [below=-.05, text style]
			{
} (0);
			
		\end{tikzpicture}
\end{minipage}
\begin{minipage}{.4\textwidth}
		\begin{tikzpicture}[auto, thick,main node/.style={circle,draw,font=\sffamily\Large\bfseries,scale=0.6},scale=0.6]
		\def\ngon{5}
		\node[regular polygon, regular polygon sides = \ngon, minimum size  = 5cm] (p) {};
		\foreach \x in {1,...,\ngon}{\node[main node] (p\x) at (p.corner \x) {$v_{\x}$};}
		\foreach \x in {1,...,\numexpr \ngon-1\relax}{
		\foreach \y in {\x,...,\ngon}{
			\draw (p\x) -- (p\y);
			}
		}
		\end{tikzpicture}
\end{minipage}
	\caption{Graphs for the {\tt flower(4,2)} strategy and {\tt completeGraph(5,1)}.}
	\label{fig:rFS}
\end{figure}
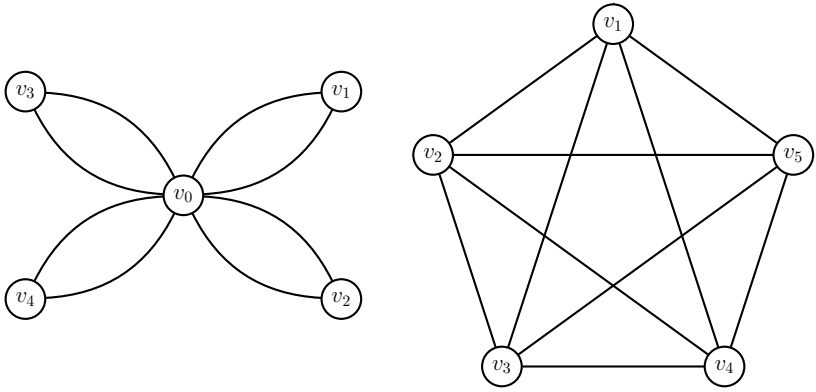

We present two graph layouts to be used for the static strategy (Figure~\ref{fig:rFS}).
\begin{description}	
\item{\verb|flower(s,t)|}
         The graph consists of a \defn{central node} $v_0$ and $s$ additional vertices (number of petals), each connected to $v_0$ by $t$ edges. 

\item{\verb|completeGraph(s,t)|}
        The graph has $s$ vertices. Every pair of vertices is connected by $t$ edges.
\end{description}

\subsubsection{Stopping criterion if a solution count is known}
\label{sec:solutioncount}
Suppose the cardinality of the fiber $\pi^{-1}(F_p)$ for a generic
value of $p$ is known.
Then a natural stopping criterion for our algorithm is to terminate when the set of known solutions $Q_i$ at any node $i$ reaches that cardinality.  
In particular, for a generic sparse system with fixed monomial support we can rely on this stopping criterion due to the BKK bound~\cite{Bernstein:BKK} that can be obtained by a \defn{mixed volume} computation.

\subsubsection{Stopping criterion if no solution count is known}
\label{sec:no-solution-count}
For a static strategy one natural stopping criterion is \defn{saturation} of the known solution correspondences along all edges. In this case, the algorithm simply can't derive any additional information. 
It also makes sense to consider a heuristic stopping criterion based
on \defn{stabilization}. The algorithm terminates when no new points are discovered in a fixed number of iterations.       
This avoids saturating correspondences unnecessarily. In particular, this could be useful if a static strategy algorithm is a part of the dynamic strategy of~\S\ref{sec:dynamic-graph-strategy}.

\begin{remark}\label{rem:trace-test}
In certain cases it is possible to provide a stopping criterion using the \defn{trace test}~\cite{SVW-trace,leykin2016trace}. This is particularly useful when there is an equation in the family $F_p(x)=0$ that describes a generic hypersurface in the parameter space, e.g., an affine linear equation with indeterminate coefficients.
In full generality, one could restrict the parameter space to a
generic line and, hence, restrict the solution variety to a
curve. Now, thinking of $F_p(x)=0$ as a system of equations
bihomogeneous in $p$ and $x$, one can use the multihomogeneous trace
test~\cite{leykin2016trace,hauenstein2015multiprojective}.  

We note that the multihomogeneous trace test complexity depends on the degree
of the solution variety, which may be significantly higher than the
degree $d$ of the covering map, where the latter is the measure of complexity for the main problem in our paper.
For instance, the system~(\ref{equation:CRN}) corresponding to the reaction network in Figure~\ref{fig:CRN} has 4~solutions, but an additional set of 11 points is necessary to execute the trace test. See \verb|example-traceCRN.m2| at~\cite{www:MonodromySolver}.     
\end{remark}

\subsubsection{Edge selection strategy}
\label{sec:edgeselection}
We propose two methods for selecting the edge $e$ in Algorithm~\ref{alg:sgs}. The default is to select an edge and direction at random.  A more
sophisticated method is to select an edge and a direction based on the potential of that selection to deliver
new information: see the discussion in Example~\ref{example:two-correspondences}.  Let $e = v_i \xleftrightarrow{(\gamma_1,\gamma_2)} v_j$ be an edge considered in the direction from $v_i$ to $v_j$.

\begin{description}
\item{$\mathtt{potentialLowerBound}$}
	equals the minimal number of new points guaranteed to be discovered by following a chosen homotopy using the maximal batch of starting points $S_i$.  That is, it equals the difference between the numbers of known unmatched points $(|Q_i|-|C_e|)-(|Q_j|-|C_e|) = |Q_i|-|Q_j|$ if this difference is positive, and $0$ otherwise.

\item{$\mathtt{potentialE}$} 
	equals the expected number of new points obtained by tracking one unmatched point along $e$. This is the ratio $d-|Q_j|\over d-|C_e|$ of undiscovered points among all unmatched points if $|Q_i|-|C_e|>0$ and $0$ otherwise. 
\end{description}

Note that $\mathtt{potentialE}$ assumes we know the cardinality of the fiber, while  $\mathtt{potentialLowerBound}$ does not depend on that piece of information.

There is a lot of freedom in choosing potentials in our algorithmic framework. The two above potentials are natural ``greedy'' choices that are easy to describe and implement. It is evident from our experiments (Tables~\ref{cyclicsevenflowergraph}~and~\ref{cyclicsevencompletegraph}) that they may order edges differently resulting in varying performance.


\subsection{An incremental dynamic graph strategy} \label{sec:dynamic-graph-strategy}
Consider a dynamic strategy that amounts to augmenting the graph once one of the above ``static'' criteria terminates Algorithm~\ref{alg:sgs} for the current graph. 
One simple way to design a dynamic stopping criterion, we call it \defn{dynamic stabilization}, is to decide how augmentation is done and fix the number of augmentation steps that the algorithm is allowed to make without increasing the solution count. 
A dynamic strategy, which is simple to implement, is one that starts with a small graph $G$ and augments it if necessary. 
\begin{algorithm}[Dynamic graph strategy] \label{alg:dgs} Let us make the same assumptions as in Algorithm~\ref{alg:sgs}.

\noindent
$(j,Q_j) = {\tt dynamicMonodromySolve}(G,x_1,{\tt stop},{\tt augment})$
\begin{algorithmic}[0]
\renewcommand{\algorithmicrequire}{\textbf{Input:}}
\renewcommand{\algorithmicensure}{\textbf{Output:}}
\REQUIRE
~\\
\begin{itemize}
\setlength{\itemsep}{0pt}
\item A graph $G$ as in Algorithm~\ref{alg:sgs}.
\item One seed solution $x_1\in\pi^{-1}(\varphi(p_1))$.
\item A stopping criterion ${\tt stop}$.
\item An augmenting procedure ${\tt augment}$.
\end{itemize}
\ENSURE A vertex $j$ in $G$ and a subset $Q_j$ of the fiber $\pi^{-1}(F_{p_j})$. 

\smallskip\hrule\smallskip

\STATE $Q_1:=\{x_1\}$ and $Q_i=\emptyset$ for $i\in 2,\dots,|V(G)|$.
\LOOP
  \STATE $(j,Q_j) = {\tt monodromySolve}(G,Q,{\tt stop})$ \COMMENT{here $Q_i$ are modified in-place and passed to the next iteration}
  \IF{{\tt stop} (i.e., stopping criterion is satisfied)}
    \RETURN $(j,Q_j)$
  \ENDIF
  \STATE $G := {\tt augment}(G)$
\ENDLOOP
\smallskip\hrule\smallskip
\end{algorithmic}
\end{algorithm}

We emphasize that the criteria described in this subsection and parts
of \S\ref{sec:no-solution-count} are \emph{heuristic} and there is a
lot of freedom in designing such. In~\S\ref{subsec:CRN} we
successfully experiment using a static stabilization criterion with some examples, for which the solution count is generally not known.

\section{Statistical analysis}
\label{sec:random_loops}
The directed cycles in the graph $G$ starting and ending at a vertex $v_1$ give elements of the fundamental group $\pi_1(B\setminus D)$, which correspond to the elements of the monodromy subgroup $M(G)$ of the monodromy group $M(\pi_1(B\setminus D))$. The latter is a subgroup of $\Sym_d$, where $d = |\pi^{-1}(F_{p_1})|$.
For example, if $G={\tt completeGraph}(2,j+1)$, then the $j$ cycles produced by edges $e_1$ and $e_2$, $e_2$ and $e_3$,..., $e_j$ and $e_{j+1}$ suffice to generate $M(G)$. The minimal number $j$ of cycles necessary to generate $M(G)$ in the general case is $\beta_1(G)$, the first Betti number of $G$ as a topological space.

For the purpose of simplifying statistical analysis, we assume that 
picking a random decorated graph $G$ with $j=\beta_1(G)$ induces uniformly and independently distributed permutations $\sigma_1, \ldots , \sigma_j \in \Sym_d$\,, where $\Sym_d$ is the symmetric group acting on the fiber $\pi^{-1}(F_{p_1})$. It would be hard in practice to achieve uniformity even when the monodromy group is a full symmetric group: see \S\ref{sec:randomization}.

\subsection{The probability of a transitive action}
Suppose the number of solutions $d$ is known and ${\tt stop}(d)$ denotes the corresponding stopping criterion. 
Our aim is to analyze the probability of producing the full solution set via Algorithm~\ref{alg:sgs} or, equivalently, the probability of
$$
{\tt dynamicMonodromySolve}(G,x_1,{\tt stop(d)},{\tt augment})
$$
terminating after at most $j$ iterations, assuming that
$\beta_1(G)=j$ at the $j$-th iteration.
This equals the probability of $\langle\sigma_1, \ldots , \sigma_j\rangle$ acting transitively, i.e., $\Pr[X_d\leq j]$
where $X_d$ is the random variable
\[
X_d = \displaystyle\inf \{ i \in \mathbb{N} \mid \langle \sigma_1, \ldots , \sigma_i \rangle \textrm{ is transitive} \}.
\]
When $d>1$ we have $\Pr [X_d =0]=0$, while $\Pr [X_d =1 ]$ is proportional to the number of $d$-cycles in the monodromy group. 
When the monodromy group is full symmetric, we can compute and give asymptotic estimates for the distribution of $X_d.$ The following theorem is a generalization of a result by Dixon, regarding the case $j=2$. The proof we give in \S\ref{sec:proof-dixon} follows the strategy of~\cite{dixon1969probability}. 
\begin{theorem}
\label{thm:dixon}
For $j\ge 2,$ $\Pr \left[ X_d \le j \right] =  1 - d^{1-j} + R_j (d),$ where the error term $R_j$ satisfies $|R_j(d)|=O(d^{-j})$.
\end{theorem}
\begin{remark}\label{rem:expected-Xd}
As a corollary, one can deduce that the expected value of $X_d$ is asymptotically finite and $\mathrm{E}[X_d]\to 2$ as $d\to\infty$.  The numerical approximations in Table~\ref{tab:tvalues} show that $\mathrm{E}[X_d]\leq 2.1033$ for all $d$. Moreover, the proof in~\S\ref{sec:proof-dixon} implies that $|R_j(d)|<C\left({d\over 2}\right)^{-j}$ with the constant $C$ not depending on $j$. Therefore, $\Pr \left[ X_d > j \right]$ decays exponentially with $j$.

Under the idealistic assumption that new cycles in the graph lead to independently and uniformly distrubuted permutations of the fiber $Q_1$, the expected Betti number needed for completion in Algorithm~\ref{alg:dgs} is at most $2.1033$. If we assume that \texttt{augment}
increases the Betti number by one by adding at most a fixed number of edges,
then \emph{the expected number of tracked paths is linear in $d$}.
\end{remark}

\begin{remark}
 We point out that Babai~\cite{babai1989probability}  proved Dixon's conjecture stating that the subgroup of $\Sym_d$ generated by two random permutations is $\Sym_d$ or $\Alt_d$ with probability $1 - d^{-1} + O(d^{-2})$. 
This shows that other subgroups are rare. However, it is easy to construct families with a transitive monodromy group that is neither full symmetric nor alternating. For example, take $x_1^2-c_1=x_2^2-c_2=0$ with irreducible solution variety and 4 solutions for generic choices of $c_1$ and $c_2$. Tracking two solutions with the same $x_1$ coordinate,  as $c_1$ and $c_2$ vary, the moving points on the tracked paths will continue to have equal projections to $x_1$. The monodromy group is $\ZZ_2\times\ZZ_2$.
\end{remark}

We reiterate that generators for the monodromy group are seldomly known \emph{a priori}. Computing them is likely to be prohibitively expensive, and the probability distribution with which our algorithm picks elements of the monodromy group is unknown, as it is prohibitively hard to analyze.


\subsection{Proof of a generalization of Dixon's theorem}
\label{sec:proof-dixon}
Fix any integer $j \ge 2.$ We wish to prove Theorem~\ref{thm:dixon} by estimating the quantity
\[
t_d = \Pr \left[ \langle \sigma_1 , \sigma_2, \ldots , \sigma_j \rangle \text{ is transitive}\right],
\]
where $\sigma_1, \ldots , \sigma_j$ are independent and uniformly distributed on $\Sym_d$.
Suppose we partition the set $\{ 1, 2, \ldots , d\}$ in such a way that there are $k_i$ classes of size $i$ for each $1\le i \le d.$ All such partitioning schemes are indexed by the set
\[
K_d = \left\{ \vec{k} \in \NN^d \mid \sum i \, k_i = d \right\} .
\]
The number of partitions corresponding to each $\vec{k} \in K_d$ is ${d!}/(\prod_{i=1}^d (i!)^{k_i} \cdot k_i!)$.
For each $\vec{k} \in K_d,$ the partition given by the orbits of $\langle \sigma_1, \ldots, \sigma_j \rangle $ is $\vec{k}$-indexed precisely when this group acts transitively on all classes of some partition associated to $\vec{k}.$ The number of tuples in $S_i^j$ with coordinates generating a group acting transitively on $\{1,\dots,i\}$ is $t_i \,(i!)^j$. Thus, we may count the set $\underbrace{\Sym_d \times \cdots \times \Sym_d}_{j}$ as
\begin{align*}
(d!)^j &= \displaystyle\sum_{\vec{k} \in K_d}  \displaystyle\frac{d!}{\prod_{i=1}^d (i!)^{k_i} \cdot k_i!} \cdot \displaystyle\prod_{i=1}^d \left( t_i \, (i!)^j  \right)^{k_i} \\
&= d! \cdot \displaystyle\sum_{\vec{k} \in K_d}\displaystyle\prod_{i=1}^d \displaystyle\frac{\left( t_i \, (i!)^{j-1} \right)^{k_i}}{k_i!}. 
\end{align*}
Let $\widehat{F}$ denote the generating function of the sequence $F(d) = (d!)^{j-1}.$ Note the formal identity
\[
\exp \left( \displaystyle\sum_{i=0}^{\infty } y_i x^i \right) =  \displaystyle\sum_{d=0}^{\infty} x^d \,  \sum_{\vec{k}  \in K_d } \prod_{i=1}^d \frac{y_i^{k_i}}{k_i!},
\]
which follows by letting  $g(x)$ denote the right hand side as a formal power series in $x,$ $f(x) = \sum y_i x^i,$ and noting the equivalent form $f' g = g'$ with $f(0)=y_0.$ We have
\begin{align*}
\displaystyle\sum_{d=1}^{\infty} d\cdot (d!)^{j-1} \, x^{d-1} &= \frac{d}{d x} \, \widehat{F} (x)  \\
&= \frac{d}{dx} \, \exp \left( \displaystyle\sum_{i=1}^{\infty } t_i \, (i!)^{j-1} x^i \right)\\
&= \Big( \displaystyle\sum_{d=0}^{\infty} (d!)^{j-1} \, x^d \Big)  \cdot \displaystyle\sum_{i=1}^{\infty } i \cdot t_i \, (i!)^{j-1} x^{i-1} \\
&= \displaystyle\sum_{d'=1}^{\infty } x^{d'-1} \, \Big( \displaystyle\sum_{i=1}^{d'}  i \cdot t_i \, \big( i! \cdot (d'-i)! \big)^{j-1} \Big)
\end{align*}
where the first equation follows by formal differentiation of the power series $\widehat{F}$, the second from the two identities above with properly substituted values for $y_i$, the third by applying the chain rule and the definition of $\widehat{F}$, and the fourth by rearranging terms by index substitution $d'=i+d$.
Upon equating coefficients of $x^{d-1}$ for $d=1,\dots$ we obtain
\begin{equation}
\label{dee_too}
d = \displaystyle\sum_{i=1}^{d} \binom{d}{i}^{1-j} \, i\,  t_i.
\end{equation}

\begin{remark}
\label{remark:computingprobabilities}
Equation (\ref{dee_too}) gives a list of linear equations in the probabilities $t_1,t_2,\dots$ allowing us to successively determine these values by backward substitution. In Table~\ref{tab:tvalues} we list some solutions for $j=2,3,4$.
\end{remark}

\begin{table}
\begin{center}
\begin{tabular}[ht]{l||l|l|l|l|}
$d$ &   $j=2$ & $j=3$ & $j=4$ & $\mathrm{E}[X_d]$\\
\hline
1  &      1
 &      1
 &      1
&0\\
2  &       0.75
&      0.875
 &      0.9375
&2\\
3  &       0.72222222
&      0.89814815
 &      0.96450617
&2.10000000\\
4  &      0.73958333
 &      0.93012153
 &      0.98262080
&2.10329381\\
5  &       0.76833333
&      0.95334722
 &      0.99115752
&2.08926525\\
10 &       0.88180398
&      0.98954768
 &      0.99898972
&2.02976996\\
20 &       0.94674288
&      0.99747856
 &      0.99987487
&2.00591026\\
30 &      0.96536852
 &      0.99888488
 &      0.99996295
& 2.00245160
\end{tabular}
\end{center}
\caption{Numerical approximations of $t_d$ --- the probability of the $j$ random permutations acting transitively on a fiber of size $d$ for $j=2,3,4$. After computing these values for larger $j$, a numerical approximation of $\mathrm{E}[X_d]$ is extracted.}
\label{tab:tvalues}
\end{table}
To complete the proof of Theorem~\ref{thm:dixon}, we introduce, as in Dixon's proof, the auxiliary quantities 
\[
\begin{array}{ccc}
r_d = d\, (1 - t_d)&  \text{and} & c_d = \displaystyle\sum_{i=1}^{d-1} \binom{d}{i}^{1-j} \, i . 
\end{array}
\]
Noting that $\binom{d}{i}^{1-j}i+\binom{d}{d-i}^{1-j}(d-i)={d\over 2}\left(\binom{d}{i}^{1-j}+\binom{d}{d-i}^{1-j}\right)$, we have
\begin{align}
{c_d\over d} =\frac{1}{2} \cdot  \displaystyle\sum_{i=1}^{d-1} \binom{d}{i}^{1-j} 
&= d^{1-j} + \left[ \binom{d}{2}^{1-j} + \frac{1}{2} \displaystyle\sum_{i=3}^{d-3} \binom{d}{i}^{1-j} \right] \nonumber\\
&\leq d^{1-j} + \left[ \binom{d}{2}^{1-j} + \frac{1}{2} (d-5) \binom{d}{3}^{1-j} \right]\label{eq:usedtwice}
\end{align}
From $j\ge 2,$ it follows that the bracketed expression in (\ref{eq:usedtwice}) is $O(d^{-j}).$
Using (\ref{dee_too}),  $t_i=1-{r_i\over i}\leq 1-{r_i\over d}$ and the definition of $c_d$, we may bound $r_d$:
\begin{align*}
r_d &= d(1-t_d)=-dt_d+d=-dt_d+\displaystyle\sum_{i=1}^{d} \binom{d}{i}^{1-j} \, i\cdot t_i \tag{3}\\
&=\displaystyle\sum_{i=1}^{d-1} \binom{d}{i}^{1-j} \, i \cdot t_i  \leq\displaystyle\sum_{i=1}^{d-1} \binom{d}{i}^{1-j} \, i \cdot (1-{r_i\over d})\\
&= \sum_{i=1}^{d-1} \binom{d}{i}^{1-j}i - {1\over d}\displaystyle\sum_{i=1}^{d-1} r_i \, \binom{d}{i}^{1-j}= c_d - {1\over d}\displaystyle\sum_{i=1}^{d-1} r_i \, \binom{d}{i}^{1-j}\leq c_d. \tag{4}
\end{align*}
To bound the error term $R_j(d):=t_d - (1-d^{1-j}),$ we consider first the case where its sign is positive. Expanding $t_d=1-{r_d\over d}$ using (3) above and $i\, t_i=i-{r_i},$
\begin{align*}
t_d - (1-d^{1-j}) &= 1-\sum_{i=1}^{d-1} \binom{d}{i}^{1-j}{it_i\over d}-1+d^{1-j}=d^{1-j}-\sum_{i=1}^{d-1} \binom{d}{i}^{1-j}{it_i\over d}\\
&=d^{1-j}-\sum_{i=1}^{d-1} \binom{d}{i}^{1-j}{i\over d}+\sum_{i=1}^{d-1} \binom{d}{i}^{1-j}{r_i\over d}\\
&=\underbrace{d^{1-j}-d^{1-j}\left({1\over d}+{d-1\over d}\right)}_{=0} \, - \, \sum_{i=2}^{d-2} \binom{d}{i}^{1-j}{i\over d} \, + \, \sum_{i=1}^{d-1} \binom{d}{i}^{1-j}{r_i\over d},
\end{align*}
and we may focus on the last summation to get, for $d\geq 2,$
\begin{align*}
t_d - (1-d^{1-j}) &\leq d^{1-j}{r_{d-1}\over d-1}+\sum_{i=2}^{d-2} \binom{d}{i}^{1-j}{r_i\over d} \tag{note $r_1=0$}\\
&\leq d^{1-j}{c_{d-1}\over d-1}+\sum_{i=2}^{d-2} \binom{d}{2}^{1-j}{c_i\over d} \tag{by (4)}\\
&\leq d^{1-j}{c_{d-1}\over d-1}+\left({d^2\over 4}\right)^{1-j}\, \sum_{i=2}^{d-2} {c_i\over i}\\
&=O(d^{2-2j})+O(d^{(2-2j)+(2-j)}) \\
&= O(d^{-j}). \tag{since $j\ge 2$}
\end{align*}
The case where $R_j(d)\le 0$ may be handled similarly, using $t_d=1-{r_d\over d}$, $r_d\leq c_d$ and what we know about the content of the bracket in (\ref{eq:usedtwice}).
\begin{align*}
-t_d + (1-d^{1-j})&= -(1-{r_d \over d})+(1-d^{1-j})= {r_d \over d}-d^{1-j}\leq {c_d \over d}-d^{1-j}\\
&\leq \left[ \binom{d}{2}^{1-j} + \frac{1}{2} (d-5) \binom{d}{3}^{1-j} \right]
= O(d^{-j}).
\end{align*}


\section{Implementation}
\label{sec:implementation}
We implement the package \verb|MonodromySolver| in Macaulay2~\cite{M2www} using the functionality of the package \verb|NumericalAlgebraicGeometry|~\cite{Leykin:NAG4M2}. The source code and examples used in the experiments in the next section are available at~\cite{www:MonodromySolver}.

The main function \verb|monodromySolve| realizes Algorithms~\ref{alg:sgs}~and~\ref{alg:dgs}, see the documentation for details and many options. The tracking of homotopy paths in our experiments is performed with the native routines implemented in the kernel of Macaulay2, however,   \verb|NumericalAlgebraicGeometry| provides an ability to outsource this core task to an alternative tracker (PHCpack or Bertini). 
Main auxiliary functions---\verb|createSeedPair|, \verb|sparseSystemFamily|, \verb|sparseMonodromySolve|, and \verb|solveSystemFamily|---are there to streamline the user's experience. The last two are blackbox routines that don't assume any knowledge of the framework described in this paper.   

The overhead of managing the data structures is supposed to be negligible compared to the cost of tracking paths. However, since our implementation uses the interpreted language of Macaulay2 for other tasks, this overhead could be sizable (up to 10\% for large examples in~\S\ref{sec:experiments}).
Nevertheless, most of our experiments are focused on measuring the \emph{number of tracked paths} as a proxy for computational complexity.

\begin{remark}
This paper's discussion focuses on linear parametric systems with a nonempty dominant component.
However, the implementation works for other cases where our framework can be applied.  

For instance, if the system is linear in parameters but has no dominant component, there may still be a unique ``component of interest'' with a straightforward way to produce a seed pair. 
This is so, for instance, in the problem of finding the degree of the variety $\SO(n)$, which we use in Table~\ref{table:monodromy-vs-world}. The point $x$ is restricted to $\SO(n)$, the special orthogonal group, which is irreducible as a variety. This results in a unique ``component of interest'' in the solution variety, the one that projects onto $\SO(n)$, see~\cite{brandt2017degree} for details.

In a yet more general case of a system that is nonlinear in parameters, it is still possible to use our software.    
We outline the theoretical issues one would need to consider in~(\ref{sec:nonlinear}) of~\S\ref{sec:generalizations}.  
\end{remark}

\subsection{Randomization}\label{sec:randomization}
Throughout the paper we refer to \emph{random} choices we make, that we assume avoid various nongeneric loci. 
For implementation purposes we make simple choices. For instance, the vertices of the graph get distributed uniformly in a cube in the base space with the exception of the seeded vertex: \verb|createSeedPair| picks $(p_0,x_0)\in B\times \CC^n$ by choosing $x$ uniformly in a cube, then choosing $p_0$ uniformly in a box in the subspace $\{ p \mid  F_p(x)=0 \}$.

A choice of probability distribution on $B$ translates to some (discrete) distribution on the symmetric group $\Sym_d$. However, it is simply too hard to analyze -- there are virtually no studies in this direction. 
We make the simplest possible assumption of uniform distribution on $\Sym_d$ in order to perform the theoretical  analysis in~\S\ref{sec:random_loops} and shed some light on \emph{why} our framework works well. There is an interesting, more involved, alternative to this assumption in~\cite{galligo2011computing,galligo2012cut}, which relies on the intuition in the case $n=1$.

\subsection{Solution count}
The BKK bound, computed via mixed volume, is used as a solution count in the examples of sparse systems in~\S\ref{subsec:cyclic} and~\S\ref{subsec:Nash}. In the latter we compute mixed volume via a closed formula that involves permanents, while the former relies on general algorithms implemented in several software packages. Our current implementation uses PHCpack~\cite{V99}, which incorporates the routines of MixedVol further developed in Hom4PS-2~\cite{lee2008hom4ps}. Other alternatives are pss5~\cite{malajovich} and Gfanlib~\cite{exactmixedvolume}. While some implementations are randomized, the latter uses symbolic perturbations to achieve exactness. The computation of the mixed volume is \emph{not} a bottleneck in our algorithm. The time spent in that preprocessing stage is negligible compared to the rest of the computation.

\subsection{Certification}\label{subsec:certification}
The reader should realize that the numerical homotopy continuation we use is driven partly by heuristics. As a post-processing step, we can certify (i.e. formally prove) the completeness and correctness of the solution set to a polynomial system computed with our main method. This is possible in the scenario when
\begin{itemize}
\itemsep=0mm
\item the parameteric system is square,
\item all solutions are regular (the Jacobian of the system is invertible), and
\item the solution count is known.
\end{itemize}
We can use Smale's $\alpha$-theory \cite[\S 8]{Alpha-Theory-book} to certify an approximation to a regular solution of a square system. In a Macaulay2 package {\tt NumericalCertification}, we implement a numerical version of an \defn{$\alpha$-test} after finding an approximate solution to certify that our solution is an \defn{approximate zero} in a rigorous sense. 

One of the main functions of {\tt NumericalCertification} is \verb|certifySolutions|, which determines whether the given solution is an approximate zero of the given polynomial system. It also produces an upper bound on the distance from the approximation to the exact solution to which it is associated.
% ...


% ...
See \verb|paper-examples/example-NashCertify.m2| at~\cite{www:MonodromySolver}, which is an example of an $\alpha$-test application to the solutions of a problem described in~\S\ref{subsec:Nash}. 
In the implementation of certification, all arithmetic and linear algebra operations are done over the field of Gaussian rationals, $\QQ[\ii]/(\ii^2+1)$. To use this certification method we first convert the coefficients of the system to Gaussian rationals, then perform certification numerically.
See~\cite{Hauenstein-Sottile:alphaCertified} for a standalone software package {\tt alphaCertified} and detailed implementation notes.  


\section{Experiments}
\label{sec:experiments}
In this section we first report on experiments with our implementation
and various examples in \S\ref{subsec:sparseexperiments} and
\S\ref{subsec:CRN}. We then investigate the completion rate of
Algorithm~\ref{alg:sgs} in \S\ref{subsec:completion}. Finally we
compare against other software in \S\ref{sec:timings}.

\subsection{Sparse polynomial systems}
\label{subsec:sparseexperiments}
The example families in this subsection have the property that the
support of the equations is fixed, while the coefficients can vary freely, as long as they are generic. We run the static graph strategy Algorithm~\ref{alg:sgs} on these examples. Our timings do not include the $\alpha$-test, which was only applied in \S\ref{subsec:Nash}.
\subsubsection{Cyclic roots}\label{subsec:cyclic}
The cyclic $n$-roots polynomial system is

\begin{equation} \label{eqcyclicsys}
 \left\{
   \begin{array}{c}
     i = 1, 2, 3, 4, \ldots, n-1: 
      \displaystyle\sum_{j=0}^{n-1} ~ \prod_{k=j}^{j+i-1}
      x_{k~{\rm mod}~n}=0 \\
     x_{0}x_{1}x_{2} \cdots x_{n-1} - 1 = 0. \\
    \end{array}
 \right.
\end{equation}
This system is commonly used to benchmark polynomial system solvers. We will study the modified system with randomized coefficients and seek solutions in $(\CC\setminus\{0\})^n$. Therefore, the solution count can be computed as the mixed volume of the Newton polytopes of the left hand sides, providing a natural stopping criterion discussed in~\S\ref{sec:solutioncount}.
This bound is 924 for cyclic-7.

Tables~\ref{cyclicsevenflowergraph} and~\ref{cyclicsevencompletegraph} 
contain averages of experimental data from running twenty trials of Algorithm~\ref{alg:sgs} on cyclic-7. 
The main measurement reported is the average number of paths tracked, as the unit of work for our algorithm is tracking a single homotopy path.
The experiments were performed with 10 different graph layouts and 3 edge selection strategies.

\begin{table}[ht]
\centering
\resizebox{\columnwidth}{!}{
\begin{tabular}{l||c|c|c|c|c|}
(\#vertices-1, edge multiplicity) & (3,2) & (4,2) & (5,2) & (3,3) & (4,3) \\ \hline
$|E(G)|$ & 6 & 8 & 10 & 9 & 12\\
$\beta_1(G)$ & 3 & 4 & 5 & 6 & 8\\
$|E(G)|\cdot 924$ &5544&7392&9240& 8316&11088 \\
\hline
completion rate     & 100\% & 100\% & 100\% & 100\% & 100\% \\ \hline
Random Edge              & 5119 & 6341 & 7544 & 6100 & 7067\\
{\tt potentialLowerBound} & 5252 & 6738 & 8086 & 6242 & 7886\\
{\tt potentialE}         & 4551 & 5626 & 6355 & 4698 & 5674\\
\end{tabular}
}
\caption{Cyclic-7 experimental results for the {\tt flower} strategy.}
\label{cyclicsevenflowergraph}
\end{table}

\begin{table}[ht]
\centering
\resizebox{\columnwidth}{!}{
\begin{tabular}{l||c|c|c|c|c|}
(\#vertices, edge multiplicity) & (2,3) & (2,4) & (2,5) & (3,2) & (4,1) \\ \hline
$|E(G)|$ & 3 & 4 & 5 & 6 & 6\\
$\beta_1(G)$ & 2 & 3 &4 & 4 & 3\\
$|E(G)|\cdot 924$ &2772&3698&4620& 5544&5544 \\
\hline
completion rate     &  65\% & 80\%  & 90\%  & 100\% & 100\% \\ \hline
Random Edge                & 2728 & 3296 & 3947 & 4805 & 5165 \\ 
{\tt potentialLowerBound}   & 2727 & 3394 & 3821 & 4688 & 5140  \\
{\tt potentialE}           & 2692 & 2964 & 2957 & 3886 & 4380 \\
\end{tabular}
}
\caption{Cyclic-7 experimental results for the {\tt completeGraph} strategy.}
\label{cyclicsevencompletegraph}
\end{table}

With respect to number of paths tracked, we see that it is an
advantage to keep the Betti number high and edge number low. 

\begin{remark}\label{rem:completion-experiments}
Computing the expected success rates ($>99\%$) using
Remark~\ref{remark:computingprobabilities}, we conclude that
the resulting permutations do not conform to the model of picking
uniformly from $S_{924}$. The completion rate depends on the choice of
strategy (compare Table~\ref{cyclicsevenflowergraph} to
Table~\ref{cyclicsevencompletegraph}). 
Nevertheless, both in theory
(assuming uniform distribution as in~\ref{sec:random_loops}) 
and
in practice (with distribution unknown to us), the completion rate does converge to~$100\%$ rapidly as
the Betti number grows. 
\end{remark}

\subsubsection{Nash equilibria}\label{subsec:Nash}
Semi-mixed multihomogeneous systems arise when one is looking for all totally
mixed Nash equilibria (TMNE) in game theory.  A specialization of mixed
volume using matrix permanents gives a concise formula for a root
count for systems arising from TMNE problems
~\cite{Emiris:2014:RCS:2608628.2608679}.  
We provide an overview of how such systems are constructed based on
~\cite{Emiris:2014:RCS:2608628.2608679}. Suppose there are $N$ players
with $m$ options each. For player $i \in \{1,\ldots,
N\}$ using option $j \in \{1,\ldots, m\}$ we have the equation $P_j^{(i)} = 0$, where

\begin{equation}
\begin{aligned}
\label{equation:NashSys}
P_j^{(i)} = \sum_{\substack{k_1,\ldots,k_{i-1},\\ k_{i+1},\ldots, k_N}}
a^{(i)}_{k_1,\ldots, k_{i-1}, j, k_{i+1}, \ldots,
  k_N}p_{k_1}^{(1)}p_{k_2}^{(2)}\cdots
p_{k_{i-1}}^{(i-1)}p_{k_{i+1}}^{(i+1)} \cdots p_{k_N}^{(N)}.
\end{aligned}
\end{equation}
The parameters $a^{(i)}_{k_1,k_2,\ldots,k_N}$ are the payoff rates for
player $i$ when players $1,\ldots,i-1,i+1,\ldots, N$ are using
options $k_1, \ldots,k_{i-1},k_{i+1},\ldots,k_N$,
respectively. Here the unknowns are $p^{(i)}_{k_j}$, representing
the probability that player $i$ will use option $k_j \in \{1,\ldots,
m\}$. There is one constraint on the probabilities for each player $i
\in \{1,\ldots, N\}$,
namely the condition that 

\begin{equation}
\label{equation:prob}
p_1^{(i)}+p_2^{(i)}+\cdots+p_m^{(i)}=1. 
\end{equation}
The system (\ref{equation:NashSys}) consists of $N\cdot m$ equations in
$N\cdot m$ unknowns. Using condition (\ref{equation:prob}) reduces the
number of unknowns to $N(m-1)$. Lastly, we eliminate the $P_j^{(i)}$
by constructing 

\begin{equation}
\label{equation:TMNESys}
P_1^{(i)}=P_2^{(i)}, ~ P_1^{(i)}=P_3^{(i)}, ~ \ldots ~ ,
P_1^{(i)}=P_m^{(i)},\quad \text{ for each } i \in \{1,\ldots, N\}.
\end{equation}
The final system is a square system of $N(m-1)$ equations in
$N(m-1)$ unknowns.

For one of our examples (\verb|paper-examples/example-Nash.m2| at \cite{www:MonodromySolver}), 
we chose the generic system of this form for $N=3$ players with
$m=3$ options for each. The result is a system of six
equations in six unknowns and 81 parameters with 10 solutions.
We also use this example to demonstrate that these solutions can be certified using {\tt NumericalCertification} (\S\ref{subsec:certification}).

\subsection{Chemical reaction networks}\label{subsec:CRN}
A family of interesting examples arises from chemical reaction network
theory.  A chemical reaction network considered under the laws of
mass-action kinetics leads to a dynamical polynomial system, the
solutions of which represent all the equilibria for the given reaction
network~\cite{2015arXiv150203188G, 2015PNAS..112.2652M}. These polynomial systems are not generically sparse and we cannot easily
compute their root count. In our experiments, we used the stabilization
stopping criterion, terminating the algorithm after a fixed number of
iterations that do not deliver new points; the default is 10 fruitless iterations. 

\begin{figure}[H]
	\centering 
		\begin{tikzpicture}
		\node (a) at (0,0) {$A$};
		\node (b) at (3,0) {$2B$};
		\node (c) at (0,-1.25) {$A+C$};
		\node (d) at (3,-1.25) {$D$};
		\node (e) at (1.5, -2.75) {$B+E$};
		\path[->,font=\scriptsize,>=angle 90]
			([yshift=2.5pt]a.east) edge node[above] {$k_2$} ([yshift=2.5pt]b.west)
			([yshift=-2.5pt]b.west) edge node[below] {$k_1$} ([yshift=-2.5pt]a.east)
			([yshift=2.5pt]c.east) edge node[above] {$k_3$} ([yshift=2.5pt]d.west)
			([yshift=-2.5pt]d.west) edge node[below] {$k_4$} ([yshift=-2.5pt]c.east)
			(d.south) edge node[below] {$k_6$} (e.north east)
			(e.north west) edge node[below] {$k_5$} (c.south);

		\end{tikzpicture}
	\caption{Chemical reaction network example.}
	\label{fig:CRN}
\end{figure}
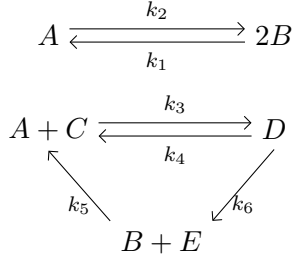

Figure~\ref{fig:CRN} gives an example of a small chemical reaction network.

Applying the laws of mass-action kinetics to the reaction network above, we obtain the polynomial
system (\ref{equation:CRN}) consisting of the corresponding steady-state and conservation
equations. Here the $k_i$'s represent the
reaction rates, $x_i$'s represent species concentrations, and the
$c_i$'s are parameters.  

\begin{equation}
\begin{aligned}
\label{equation:CRN}
\dot{x_A} & = k_1x_B^2-k_2x_A-k_3x_Ax_C+k_4x_D+k_5x_Bx_E  \\
\dot{x_B} & = 2k_1x_A-2k_2x_B^2+k_4x_D-k_5x_Bx_E \\
\dot{x_C} & = -k_3x_Ax_C+k_4x_D+k_5x_Bx_E \\
\dot{x_D} & = k_3x_Ax_C-(k_4+k_6)x_D  \\
\dot{x_E} & = -k_5x_Bx_E+k_6x_D    \\
 0 & = 2x_A+x_B-x_C+x_D-c_1 \\
0 & = -2x_A-x_B+2x_C+x_E-c_2    
\end{aligned}
\end{equation}
Typically, systems resulting from chemical reaction networks will be
overdetermined. With the current implementation one needs to either
square the system or use a homotopy tracker that supports following a
homotopy in a space of overdetermined systems.  

Although we may obtain large systems, they
typically have very low root counts compared to the sparse case.  
The polynomial system (\ref{equation:CRN}) has four solutions. 
A larger example is the \defn{wnt signaling pathway} from Systems
Biology~\cite{2015arXiv150203188G} consisting of 19 polynomial
equations with 9 solutions. All 9 solutions are obtained in less than a
second with Algorithm~\ref{alg:sgs}.

\subsection{Completion rate}\label{subsec:completion}
We investigate the completion rate of Algorithm~\ref{alg:sgs} for the Katsura family parametrized by $n$ with fixed support and generically chosen coefficients.
Tables~\ref{katsuraflowergraph} and~\ref{katsuracompletegraph} contain 
the percentage of successes from 500 runs with distinct random seeds. In Table~\ref{tab:katsuraexpected}, we show the computed expected values using Remark~\ref{remark:computingprobabilities}.

\begin{table}[ht]
\centering
\begin{tabular}{l|l||c|c|c|c|c|}
\multicolumn{2}{c||}{} &  \multicolumn{5}{c|}{(\#vertices-1, edge multiplicity)} \\ \hline
$n$     & BKK Bound & (3,2) & (4,2) & (5,2) & (3,3) & (4,3) \\ 
           & & $\beta_1=3$ & $\beta_1=4$ & $\beta_1=5$ & $\beta_1=6$& $\beta_1=8$ \\ \hline
5 & 12 & 96.4\% & 99.4\% & 99.6\% & 100\% & 99.8\% \\
6 & 30 & 98.6\% & 100\% & 99.8\% & 100\% & 99.6\% \\
7 & 54 & 97.6\% & 98.8\% & 99.4\% & 99.4\% & 98.4\% \\
8 & 126 & 99.2\% & 99.8\% & 99.6\% & 99.8\% & 99.8\% \\
9 & 240 & 98.8\% & 99.6\% & 98.4\% & 98.4\% & 98.6\% \\
10 & 504 & 98.6\% & 98.8\% & 99.2\% & 99.4\% & 98.8\% \\
\end{tabular}
\caption{Katsura-$(n-1)$ for the {\tt flower} strategy.}
\label{katsuraflowergraph}
\end{table}

\begin{table}[ht]
\centering
\begin{tabular}{l|l||c|c|c|c|c|}
\multicolumn{2}{c||}{} &  \multicolumn{5}{c|}{(\#vertices, edge multiplicity)} \\ \hline
$n$ & BKK Bound & (2,3) & (2,4) & (2,5) & (3,2) & (4,1) \\
           && $\beta_1=2$ & $\beta_1=3$ & $\beta_1=4$ & $\beta_1=4$& $\beta_1=3$ \\ \hline
5 & 12 & 65.6\% & 88.2\% & 95\% & 99.2\% & 98\% \\
6 & 30 & 77.4\% & 95.2\% & 99\% & 99.8\% & 99.6\% \\
7 & 54 & 74.4\% & 96.2\% & 99.2\% & 99.6\% & 99.8\% \\
8 & 126 & 81.8\% & 97\% & 99.2\% & 100\% & 99.8\% \\
9 & 240 & 85.2\% & 97.6\% & 99.4\% & 99\% & 98.2\% \\
10 & 504 & 89.2\% & 98.2\% & 99.2\% & 99.4\% & 99\% \\
\end{tabular}
\caption{Katsura-$(n-1)$ for the {\tt completeGraph} strategy.}
\label{katsuracompletegraph}
\end{table}

\begin{table}[ht]
\centering
\begin{tabular}{r||r|r|r|}
 $d$    & $\beta_1=2$ & $\beta_1=3$ & $\beta_1\geq 4$ \\ \hline
12  & 90.5\% & $99.3 \%$ & 100.0\%\\
30  & 96.5\% & $99.9 \%$ & 100.0\%\\
54  & 98.1\% & 100.0\% & 100.0\%\\
126 & 99.2\% & 100.0\% & 100.0\%\\
240 & 99.6\% & 100.0\% & 100.0\%\\
504 & 99.8\% & 100.0\% & 100.0\%\\
\end{tabular}
\caption{Rounded expected probability of success assuming uniform distribuition of permutations and full monodromy group.}
\label{tab:katsuraexpected}
\end{table}

For $\beta_1\geq 3$ the observed success rates approach the expected
values of Table~\ref{tab:katsuraexpected}. We note that the {\tt
  flower} strategy is again closest to the estimates.
We do not expect the numbers produced in experiments to match the numbers in Table~\ref{tab:tvalues}, since the assumptions made for that statistical analysis are quite idealistic; however, both the analysis and experiments show that the probability of success approaches 100\% rapidly as the number of solutions grows and the first Betti number increases.

\subsection{Timings and comparison with other solvers}\label{sec:timings}
All timings appearing in this section are done on \emph{one} thread and on the \emph{same} machine.
Remarks~\ref{rem:complexity}~and~\ref{rem:expected-Xd} show that we should expect the number of tracked paths in Algorithms~\ref{alg:sgs}~and~\ref{alg:dgs} to be linear (with a small constant!) in the number of solutions of the system. 
In this section we highlight the practicality of our approach in two ways.


Firstly, 
the monodromy method dramatically extends our computational ability for systems where the solution count turns out to be significantly smaller than the count corresponding to a more general family, for example, BKK count for sparse systems. 
This means that the existing blackbox methods, whose complexity relies on a larger count, are likely to spend significantly more time in computation compared to our approach.
In Table~\ref{table:monodromy-vs-world}, we collect timings on several challenging examples mentioned in recent literature where smaller solution counts are known, thus providing us with rigorous test cases for our heuristic stopping criterion. 
The first system in the table is that of the wnt signaling pathway reaction network mentioned in~\S\ref{subsec:CRN}. The others come from the problem of computing the degree of $\SO(n)$, the special orthogonal group, as a variety~\cite{brandt2017degree}.

Below is a list of comments on the setup:
\begin{itemize}
\item For our implementation we chose small graphs with $\beta_1\leq 4$ and the random edge selection strategy. The stopping criterion is ``stabilization'' as discussed in~\S\ref{sec:no-solution-count}. 

\item While the blackbox solver of PHCpack ultimately performs polyhedral homotopy continuation, Bertini relies by default on an equation-by-equation technique dubbed \defn{regeneration} (see~\cite{Bertini-book}). The latter may be faster than the former in certain cases, which this series of examples shows.
\end{itemize}

\begin{table}[h]
\centering
\begin{tabular}{l||c|c|c|c|c|}
problem     & wnt & $\SO(4)$ & $\SO(5)$  & $\SO(6)$  & $\SO(7)$ \\
\hline
count       & 9           &  40      & 384       & 4768      & 111616 \\
\hline
MonodromySolver & 0.52    & 4        & 23        & 528       & 42791   \\
Bertini         & 42      & 81       & 10605     & out of memory &     \\
PHCpack         & 862     & 103      & $>$ one day &  &         \\
\end{tabular}
\caption{Examples with solution count smaller than BKK bound (timings in seconds).}
\label{table:monodromy-vs-world}
\end{table}

Secondly, when the solution count is given by the BKK bound our method is a viable alternative to polyhedral homotopy solvers, since the number of paths we track is linear in the number of solutions. 
The timings on a few large benchmark problems of our current implementation and several other software packages are in Table~\ref{table:timings}. 
Our goal in the rest of this section is to show that our running times are in the same ballpark as polyhedral homotopies. 
 
Below is a list of comments on the setup:
\begin{itemize}
\item For our implementation we chose two small graphs and default (random) edge selection strategy.

\item For PHCpack there is a way to launch a mixed volume computation with the option of creating a system with the same support and random coefficients together with its solutions. This is the option we are using; the blackbox computation takes a little longer.

\item HOM4PS2~\cite{lee2008hom4ps} is not open source unlike all other software mentioned here. 
(We use HOM4PS2 stock examples for all systems and call its blackbox polyhedral homotopies solver.)
HOM4PS2 may use just-in-time compilation of straight-line programs used for evaluation, which speeds up computations considerably. (PHCpack does not use this technique; neither does our software, but our preliminary experiments in Macaulay2 show a potential for a 10- to 20-fold speed up over our currently reported timings.)
\end{itemize}

\begin{table}[h]
\centering
\begin{tabular}{l||c|c|c|}
problem                          & cyclic-10 & cyclic-11 & noon-10\\
\hline
BKK bound                       & 35940     &  184756   & 59029 \\
\hline
{\tt completeGraph}(2,3) & 610 & 7747& failed \\
              & {\small(107820 paths)} &  {\small(540155 paths)} & {\small (59001 solutions)}\\
{\tt completeGraph}(2,4) & 740 & 8450  & 935 \\
            &{\small(129910 paths)} &{\small(737432 paths)} & {\small(236051 paths)}\\
PHCpack                          & 538         &4256                                        & 751 \\
HOM4PS2                         & 62           &410                                         & 120\\
\end{tabular}
\caption{Software timings on large examples (in seconds).}
\label{table:timings}
\end{table}

\begin{remark} For large examples, assuming the probabilistic model leading to Theorem~\ref{thm:dixon} and Remark~\ref{rem:expected-Xd}, the probability of success should be extremely close to 100\% even for a random graph with $\beta_1 = 2$. The run of noon-10, which is an example of neural network model from~\cite{noonburg1989neural}, demostrates an unlikely but possible failure for $\beta_1 = 2$ followed by success at $\beta_1 = 3$.      
\end{remark}

On the examples in Table~\ref{table:timings}, we also ran the blackbox solvers 
of Bertini and NumericalAlgebraicGeometry~\cite{Leykin:NAG4M2}, which use the total-degree 
homotopy. Both were able to finish noon-10 with timings similar to the table,
but all other problems took longer than a day.
This is expected, as the BKK bound of noon-10 is only slightly	
sharper than the B\'ezout bound.

\begin{remark} In comparison with the naive dynamic strategy (\S\ref{subsec:naive}) our framework loses slightly only in one aspect: memory consumption. For a problem with $d$ solutions the naive approach stores up to (and typically close to) $2d$ points. The number of points our approach stores is up to (and typically considerably fewer than) $d$ times the number of vertices. For instance,  it is up to $4d$ points in all runs in Table~\ref{table:timings}.

The number of tracked paths is significantly lower in our framework: for example, the naive strategy tracks about 7500 paths on average for cyclic-7.  Even before looking at Table~\ref{cyclicsevenflowergraph} it is clear that running the \verb|flower| strategy in combination with the incremental dynamic strategy of~\S\ref{sec:dynamic-graph-strategy} guarantees to dominate the naive strategy.
\end{remark}
\section{Generalizations}\label{sec:generalizations}
While we propose a more general algorithmic framework, a concurrent goal of this paper is to demonstrate that significant practical advantages are already apparent when we apply a relatively simple implementation and analysis to \emph{simple problems} (linearly parametrized families). The following topics thus lie outside the scope of this article, but seem deserving of further study:
\begin{enumerate}

\item \label{sec:failures}
One advantage of the MS approach is that it can tolerate numerical failures of the underlying homotopy tracker. In fact, we already implemented a simple failure resistant mechanism and it successfully tolerates a few failures that arise in some runs for large test examples in~\S\ref{sec:timings}. A natural extension of this paper's statistical analysis would be to model the algorithm's performance in the presence of failures. 

\item Ideally, heuristics such as edge potentials should incorporate information such as the failures discussed above. It is also of interest to adapt potentials to the \emph{parallel} setting discussed below.

\item The parallelization of the MS approach is not as straightforward as that of other homotopy continuation methods. The question of when speedups close to linear can be achieved should be addressed. 

\item \label{sec:nonlinear} Consider the generalized setup in which the base space $B$ is an irreducible variety and the family is given by a rational map from $P$ into a space of systems. To apply our general framework, a major requirement is to find an effective way to parametrize a curve between two points of $P.$ This parametrization would conceivably depend on the nature of the problem being considered. Certain other ingredients are also likely to be problem-specific---for instance, even in the case of $P = \CC^m,$ the construction of the initial seed $(p_0,x_0)$ is complicated by the possibility that the systems' coefficients are nonlinear in the parameters. Nonetheless, this is one of the strengths of the MS framework---once all required ``oracles'' are supplied, the procedures become effective.

\item \label{sec:Galois} In the classical language of enumerative geometry, the monodromy groups we consider are isomorphic to Galois groups of incidence varieties (essentially solution varieties in our terminology). For a large class of Schubert problems and other interesting incidence varieties, the associated Galois group turns out to be the full symmetric group.~\cite{Leykin-Sottile:HoG} A suitable modification of our dynamic strategy is one practical approach to verifying this in conjectural cases. 

\item Our paper demonstrates the strength of our method relative to other techniques such as polyhedral homotopy and regeneration. Building on our framework, one could use polyhedral homotopy as a subroutine to quickly populate a partial solution set (quickly discarding any path that becomes poorly conditioned). Further advantages may be achievable by using different techniqes in parallel. These and other hybrid approaches have the potential to produce even faster and more robust blackbox solvers.
\end{enumerate}


\bibliography{refs}

\begin{thebibliography}{10}

\bibitem{babai1989probability}
L{\'a}szl{\'o} Babai.
\newblock The probability of generating the symmetric group.
\newblock {\em Journal of Combinatorial Theory, Series A}, 52(1):148--153,
  1989.

\bibitem{Bertini-book}
Daniel~J. Bates, Jonathan~D. Hauenstein, Andrew~J. Sommese, and Charles~W.
  Wampler.
\newblock {\em Numerically solving polynomial systems with {B}ertini},
  volume~25.
\newblock SIAM, 2013.

\bibitem{Bernstein:BKK}
D.~N. Bernstein.
\newblock The number of roots of a system of equations.
\newblock {\em Funkcional. Anal. i Prilo\v zen}, 9(3):1--4, 1975.

\bibitem{Alpha-Theory-book}
Lenore Blum, Felipe Cucker, Michael Shub, and Steve Smale.
\newblock {\em Complexity and real computation}.
\newblock Springer-Verlag, New York, 1998.

\bibitem{brandt2017degree}
M.~Brandt, J.~Bruce, T.~Brysiewicz, R.~Krone, and E.~Robeva.
\newblock The degree of $\text{SO}(n,\mathbb{C})$.
\newblock In G.~Smith and B.~Sturmfels, editors, {\em Combinatorial Algebraic
  Geometry}, volume~80 of {\em Fields Institute Communications}, pages
  229--246. Springer, 2017.

\bibitem{Chen-Kileel:NumericalImplicitization}
Justin Chen and Joe Kileel.
\newblock {Numerical Implicitization for Macaulay2}.
\newblock {\em arXiv preprint arXiv:1610.03034}, 2016.

\bibitem{del2015critical}
Abraham~M. del Campo and Jose~I. Rodriguez.
\newblock Critical points via monodromy and local methods.
\newblock {\em Journal of Symbolic Computation}, 79:559--574, 2017.

\bibitem{dixon1969probability}
John~D. Dixon.
\newblock The probability of generating the symmetric group.
\newblock {\em Mathematische Zeitschrift}, 110(3):199--205, 1969.

\bibitem{www:MonodromySolver}
T.~Duff, C.~Hill, A.~Jensen, K.~Lee, A.~Leykin, and J.~Sommars.
\newblock {MonodromySolver:} a {Macaulay2} package for solving polynomial
  systems via homotopy continuation and monodromy.
\newblock Available at
  http://people.math.gatech.edu/$\sim$aleykin3/MonodromySolver.

\bibitem{Emiris:2014:RCS:2608628.2608679}
Ioannis~Z. Emiris and Raimundas Vidunas.
\newblock Root counts of semi-mixed systems, and an application to counting
  nash equilibria.
\newblock In {\em Proceedings of the 39th International Symposium on Symbolic
  and Algebraic Computation}, ISSAC '14, pages 154--161, New York, NY, USA,
  2014. ACM.

\bibitem{galligo2012cut}
Andr{\'e} Galligo and Laurent Miclo.
\newblock On the cut-off phenomenon for the transitivity of randomly generated
  subgroups.
\newblock {\em Random Structures \& Algorithms}, 40(2):182--219, 2012.

\bibitem{galligo2011computing}
Andr{\'e} Galligo and Adrien Poteaux.
\newblock Computing monodromy via continuation methods on random {R}iemann
  surfaces.
\newblock {\em Theoretical Computer Science}, 412(16):1492--1507, 2011.

\bibitem{M2www}
Daniel~R. Grayson and Michael~E. Stillman.
\newblock Macaulay2, a software system for research in algebraic geometry.
\newblock Available at http://www.math.uiuc.edu/Macaulay2/.

\bibitem{2015arXiv150203188G}
E.~{Gross}, H.~A. {Harrington}, Z.~{Rosen}, and B.~{Sturmfels}.
\newblock Algebraic systems biology: A case study for the {W}nt pathway.
\newblock {\em Bulletin of Mathematical Biology}, 78(1):21--51, 2016.

\bibitem{hauenstein2015multiprojective}
Jonathan~D. Hauenstein and Jose~I. Rodriguez.
\newblock Multiprojective witness sets and a trace test.
\newblock {\em arXiv preprint arXiv:1507.07069}, 2015.

\bibitem{hauenstein2016numerical}
Jonathan~D. Hauenstein, Jose~I. Rodriguez, and Frank Sottile.
\newblock Numerical computation of {G}alois groups.
\newblock {\em Found Comput Math}, 2017.

\bibitem{Hauenstein-Sottile:alphaCertified}
Jonathan~D. Hauenstein and Frank Sottile.
\newblock {Algorithm 921: alphaCertified: Certifying} solutions to polynomial
  systems.
\newblock {\em ACM Trans. Math. Softw.}, 38(4):28:1--28:20, August 2012.

\bibitem{HuberSturmfels:PolyhedralHomotopies}
Birkett Huber and Bernd Sturmfels.
\newblock A polyhedral method for solving sparse polynomial systems.
\newblock {\em Math. Comp.}, 64(212):1541--1555, 1995.

\bibitem{exactmixedvolume}
Anders~N. Jensen.
\newblock An implementation of exact mixed volume computation.
\newblock In {\em Mathematical Software - {ICMS} 2016 - 5th International
  Conference, Berlin, Germany, July 11-14, 2016, Proceedings}, pages 198--205,
  2016.

\bibitem{lee2008hom4ps}
Tsung-Lin Lee, Tien-Yien Li, and Chih-Hsiung Tsai.
\newblock {HOM4PS-2.0}: A software package for solving polynomial systems by
  the polyhedral homotopy continuation method.
\newblock {\em Computing}, 83(2-3):109--133, 2008.

\bibitem{Leykin:NAG4M2}
Anton Leykin.
\newblock Numerical algebraic geometry.
\newblock {\em The Journal of Software for Algebra and Geometry}, 3:5--10,
  2011.

\bibitem{leykin2016trace}
Anton Leykin, Jose~I. Rodriguez, and Frank Sottile.
\newblock Trace test.
\newblock {\em arXiv preprint arXiv:1608.00540}, 2016.

\bibitem{Leykin-Sottile:HoG}
Anton Leykin and Frank Sottile.
\newblock Galois groups of {S}chubert problems via homotopy computation.
\newblock {\em Math. Comp.}, 78(267):1749--1765, 2009.

\bibitem{Ley-Ver-new-monodromy-05}
Anton Leykin and Jan Verschelde.
\newblock Decomposing solution sets of polynomial systems: A new parallel
  monodromy breakup algorithm.
\newblock {\em International Journal of Computational Science and Engineering},
  4(2):94--101, 2009.

\bibitem{2015PNAS..112.2652M}
A.~L. {MacLean}, Z.~{Rosen}, H.~M. {Byrne}, and H.~A. {Harrington}.
\newblock {Parameter-free methods distinguish Wnt pathway models and guide
  design of experiment}.
\newblock {\em Proceedings of the National Academy of Science}, 112:2652--2657,
  March 2015.

\bibitem{malajovich}
Gregorio Malajovich.
\newblock Computing mixed volume and all mixed cells in quermassintegral time.
\newblock {\em Found. Comput. Math.}, 17(5):1293--1334, 2017.

\bibitem{Morgan87}
Alexander Morgan.
\newblock {\em Solving polynomial systems using continuation for engineering
  and scientific problems}.
\newblock Prentice Hall Inc., Englewood Cliffs, NJ, 1987.

\bibitem{noonburg1989neural}
V.~W. Noonburg.
\newblock A neural network modeled by an adaptive {L}otka-{V}olterra system.
\newblock {\em SIAM Journal on Applied Mathematics}, 49(6):1779--1792, 1989.

\bibitem{SVW2001:monodromy}
A.~J. Sommese, J.~Verschelde, and C.~W. Wampler.
\newblock {\em Using monodromy to decompose solution sets of polynomial systems
  into irreducible components}, pages 297--315.
\newblock Springer Netherlands, Dordrecht, 2001.

\bibitem{SVW-trace}
A.~J. Sommese, J.~Verschelde, and C.~W. Wampler.
\newblock Symmetric functions applied to decomposing solution sets of
  polynomial systems.
\newblock {\em SIAM J.\ Numer.\ Anal.}, 40(6):2026--2046, 2002.

\bibitem{SVW9}
Andrew~J. Sommese, Jan Verschelde, and Charles~W. Wampler.
\newblock Introduction to numerical algebraic geometry.
\newblock In Alicia Dickenstein and Ioannis~Z. Emiris, editors, {\em Solving
  Polynomial Equations: Foundations, Algorithms, and Applications}, pages
  301--337. Springer Berlin Heidelberg, Berlin, Heidelberg, 2005.

\bibitem{Sommese-Wampler-book-05}
Andrew~J. Sommese and Charles~W. Wampler, II.
\newblock {\em The numerical solution of systems of polynomials}.
\newblock World Scientific Publishing Co. Pte. Ltd., Hackensack, NJ, 2005.

\bibitem{V99}
Jan Verschelde.
\newblock Algorithm 795: {PHC}pack: A general-purpose solver for polynomial
  systems by homotopy continuation.
\newblock {\em ACM Trans. Math. Softw.}, 25(2):251--276, 1999.

\bibitem{Verschelde-Verlinden-Cools}
Jan Verschelde, Pierre Verlinden, and Ronald Cools.
\newblock Homotopies exploiting {N}ewton polytopes for solving sparse
  polynomial systems.
\newblock {\em SIAM J. Numer. Anal.}, 31(3):915--930, June 1994.

\end{thebibliography}
\bibliographystyle{plain}

\end{document}